\magnification=\magstep1 
\overfullrule=0pt  \voffset 0.4in
\input amssym.def
\def\QED{\vrule height6pt width6pt depth0pt}

 \def\la{\lambda}  \def\ga{\gamma} 
\def\si{\sigma} \def\eps{\epsilon}  \def\Z{{\Bbb Z}}
\font\huge=cmr10 scaled \magstep2    \def\C{{\Bbb C}}  \def\Q{{\Bbb Q}}
\font\smcap=cmcsc10        \def\R{{\Bbb R}} 
\def\eqde{:=}   \def\io{{\iota}}   \def\crprod{{\prod}}  
\def\k{{\kappa}}  \def\i{{\rm i}}  \def\L{{\Lambda}}  \def\N{{\cal N}}
\def\crprod{{\times\!\vrule height5pt depth0pt width0.4pt\,}}
\def\vp{\varphi}  \def\E{{\cal E}}

\centerline{{\bf \huge  Modular data: the algebraic combinatorics}}
\bigskip\centerline{{\bf \huge  of conformal field theory}}
\bigskip \bigskip   \centerline{Terry Gannon}\medskip
\centerline{{\it Department of 
Mathematical Sciences, University of Alberta,}}
\centerline{{\it  Edmonton, Canada, T6G 1G8}}
\smallskip                                    
\centerline{{e-mail: tgannon@math.ualberta.ca}}
\bigskip\medskip

\noindent{{\bf Abstract.}} 
This paper is primarily intended as an introduction for mathematicians to some
of the rich algebraic combinatorics arising in for instance conformal
field theory (CFT). It tries to refine, modernise, and bridge the gap
between papers [4] and [39]. Our paper is essentially self-contained, apart from some of the
background motivation (Section I) and examples (Section III) which are included to give the reader
a sense of the context. Detailed proofs will appear elsewhere. The theory is
still a work-in-progress, and emphasis is given here to several open
questions and problems.\bigskip\medskip

\centerline{{\bf I. Introduction}}\medskip

In Segal's axioms of CFT [83], any Riemann surface with boundary is assigned
a certain linear homomorphism. Roughly speaking, Borcherds [13] and
Frenkel-Lepowsky-Meurman [37] axiomatised this data corresponding to
a sphere with 3 disks removed, and the result is called a vertex operator
algebra. Here
we do the same with the data corresponding to a torus (and to a lesser extent
a cylinder). The result is considerably simpler, as we shall see.

{\it Moonshine} in its more general sense involves the assignment of modular
(automorphic) functions or forms to certain algebraic
structures, e.g.\ theta functions to lattices, or vector-valued Jacobi forms to
affine algebras, or Hauptmoduls to the Monster. 
This paper explores an important facet of Moonshine theory:
the associated modular group representation. From this perspective,
{\it Monstrous} Moonshine [14] is maximally uninteresting: the corresponding
representation is completely trivial! 

Let's focus now on the former context. 
Do not be put-off if this
introductory section contains many unfamiliar terms. This section is
motivational, supplying some of the background physical context, and
many of the terms here will be mathematically addressed in later sections.
It is intended to be skimmed.

A rational conformal field theory
(RCFT) has two vertex operator
algebras (VOAs) ${\cal V},{\cal V}'$. For simplicity we will take them
to be isomorphic (otherwise the RCFT is called `heterotic'). The VOA
${\cal V}$ will have finitely many irreducible
modules $A$. Consider their (normalised) characters 
$${\rm ch}_A(\tau)=q^{-c/24}\,{\rm Tr}_Aq^{L_0}\eqno(1.1)$$
where $c$ is the rank of the VOA and $q=e^{2\pi\i\tau}$, for $\tau$ in the
upper half-plane ${\Bbb H}$.
A VOA ${\cal V}$ is (among other things) a vector space with a grading
given by the eigenspaces of the operator $L_0$; (1.1) defines the character
to be obtained from the induced $L_0$-grading on the ${\cal V}$-modules $A$.
These characters yield a representation of the modular group SL$_2(\Z)$ of the
torus, given by
its familiar action on ${\Bbb H}$ via fractional linear
transformations. In particular, we can define matrices $S$ and $T$ by
$${\rm ch}_A(-1/\tau)=\sum_{B}S_{AB}\,{\rm ch}_B(\tau),\qquad
{\rm ch}_A(\tau+1)=\sum_BT_{AB}\,{\rm ch}_B(\tau)\ ;\eqno(1.2a)$$
this representation sends
$$\left(\matrix{0&-1\cr 1&0}\right)\mapsto S,\qquad
\left(\matrix{1&1\cr 0&1}\right)\mapsto T\ .\eqno(1.2b)$$
We call this representation the {\it modular data} of the RCFT. It has
some interesting properties, as we shall see. For example, in
Monstrous Moonshine the relevant VOA is the Moonshine module
$V^\natural$. There is only one irreducible module of $V^\natural$,
namely itself, and its character $j(\tau)-744$ is invariant under SL$_2(\Z)$.

Incidentally, there is in
RCFT and related areas a (projective) representation of each mapping
class group --- see e.g.\ [2,43,73,85] and references therein.
These groups play the role
of modular group, for any Riemann surface. Their representations
coming from e.g.\ RCFT are still poorly understood, and certainly deserve more attention, but in this
paper we will consider only SL$_2(\Z)$ (i.e.\ the unpunctured torus).

Strictly speaking we need linear
independence of our characters, which means considering the `1-point
functions' 
$${\rm ch}_A(\tau,u)=q^{-c/24}\,{\rm Tr}_A(q^{L_0}\,o(u))$$
--- this is why SL$_2(\Z)$ and not PSL$_2(\Z)$ arises here ---
but for simplicity we will ignore this technicality in the following.

In physical parlance, the two VOAs are the (right- and left-moving)
algebras of (chiral) observables.
The observables operate on the space ${\cal H}$ of
physical states of the theory; i.e.\ ${\cal H}$ carries a representation of
${\cal V}\otimes{\cal V}$. The irreducible modules $A\otimes A'$ of
${\cal V}\otimes {\cal V}$ in ${\cal H}$ are
labelled by the {\it primary fields} --- special states
$|\phi,\phi'\rangle$ in ${\cal H}$
which play the role of highest weight vectors. More precisely, the primary
field will be a vertex
operator $Y(\phi,z)$ and the ground state $|\phi\rangle$ will be the
state created by the primary field at time $t=-\infty$: 
$|\phi\rangle={\rm lim}_{z\rightarrow 0}Y(\phi,z)|0\rangle$. The VOA ${\cal
V}$ acting on the (chiral) primary field $|\phi\rangle$ generates the module
$A=A_\phi$ (and similarly for $\phi'$). 
The characters ch$_A$ form a basis for the
vector space of  0-point 1-loop conformal blocks (see (3.7) with $g=1,t=0$).

Modular data is a fundamental ingredient of the RCFT. It appears for
instance in Verlinde's formula (2.1), which gives (by
definition) the structure constants for the fusion ring. It also constrains
the torus partition function ${\cal Z}$:
$${\cal Z}(\tau)=q^{-c/24}\,\overline{q}^{-c/24}\,{\rm Tr}_{\cal H}\,
q^{L_0}\,\overline{q}^{{L}'_0}\eqno(1.3a)$$
where $\overline{q}$ is the complex conjugate of $q$. Now as mentioned
above, ${\cal H}$ has the decomposition
$${\cal H}=\oplus_{A,B}\,M_{AB}\,A\otimes B\eqno(1.3b)$$
into ${\cal V}$-modules, where the $M_{AB}$ are multiplicities, and so
$${\cal Z}(\tau)=\sum_{A,B}M_{AB}\,{\rm ch}_A(\tau)\,\overline{{\rm
    ch}_B(\tau)}\eqno(1.3c)$$
Physically, ${\cal Z}$ is the 1-loop vacuum-to-vacuum amplitude of the
    closed string (or
rather, the  amplitude would be $\int {\cal Z}(\tau)\,d\tau$).
`Amplitudes' are the fundamental numerical quantities in quantum 
theories, from which the experimentally determinable  probabilities are obtained.
In Segal's formalism, the torus ${\Bbb C}/(\Z+\tau\Z)$ is assigned the
homomorphism ${\Bbb C}\rightarrow{\Bbb C}$ corresponding to multiplication
by ${\cal Z}(\tau)$.
We will see in Section V that  ${\cal Z}$ must be invariant under the action (1.2a)
 of the modular group SL$_2(\Z)$, and so we call it (or equivalently its matrix $M$ of multiplicities)
a {\it modular invariant}. 

Another elementary but fundamental quantity
is the 1-loop vacuum-to-vacuum amplitude ${\cal Z}_{\alpha\beta}$ of
the open string, to whose ends
are attached `boundary states'  $|\alpha\rangle,|\beta\rangle$ ---
this cylindrical partition function  looks like
$${\cal Z}_{\alpha\beta}(t)=\sum_A{\cal N}_{A\alpha}^\beta\,{\rm ch}_A(\i t)
\eqno(1.4)$$
where these multiplicities ${\cal N}_{A\alpha}^\beta$ have something to do 
with Verlinde's formula (2.1). These functions ${\cal Z}_{\alpha\beta}$ (or equivalently
their matrices $({\cal N}_A)_{\alpha\beta}={\cal N}_{A\alpha}^\beta$ of coefficients)
are called {\it fusion graphs} or {\it NIM-reps}, for reasons that
will be explained in Section V.

We define  modular invariants and NIM-reps axiomatically in Section  V.
Classifying them is
essentially the same as classifying (boundary) RCFTs, and is an
interesting and accessible challenge. All of this will be explained more thoroughly
and rigourously in the course of this paper.

In this paper we survey the basic theory and examples of modular
data and fusion rings. Then we sketch the basic theory of modular
invariants and NIM-reps. Finally, we specialise to the modular data associated to
affine Kac-Moody algebras, and discuss what is known about their
modular invariant (and NIM-rep) classifications. A familiarity with RCFT is not needed
to read this paper (apart from this introduction!).

The theory of fusion rings (and modular data) in its purest form is the study of the 
algebraic consequences of requiring structure constants to obey the 
constraints of positivity and integrality, as well as imposing some sort of
self-duality condition identifying the ring with its dual.
But one of the thoughts running through this note is that we don't know yet its
correct definition. In the next section is given the most standard definition,
but surely it can be improved. How to determine the correct definition is clear:
we probe it from the `inside' --- i.e.\ with strange examples which we probably
want to call
modular data --- and also from the `outside' --- i.e.\ with examples probably
too dangerous to include in the fold. Some of these critical examples
will be described below.

\smallskip\noindent{{\smcap Notational remarks:}} Throughout the paper we let $\Z_{\ge}$
denote the nonnegative integers, and $\overline{x}$ denote the complex
conjugate of $x$. The transpose of a matrix $A$ will be written
$A^t$.

\vfill\eject
\centerline{{\bf II. Modular data and fusion rings}}\medskip

The most basic structure considered in this paper is that of modular
data; the particular variant studied here --- and the most common one 
in the literature --- is given in Definition 1. But there are
alternatives, and a natural general one is given by {\bf MD1}$'$, {\bf
MD2}$'$, {\bf MD3}, and {\bf MD4}. In the more limited context of
e.g.\ RCFT, the appropriate axioms are {\bf MD1}, {\bf MD2}$'$, and
{\bf MD3}-{\bf MD6}.

\medskip\noindent{{\bf Definition 1.}} Let $\Phi$ be a finite set of labels, one of which --- we will denote it 0
and call it the `identity' ---
is distinguished. By {\it modular data} we mean matrices
$S=(S_{ab})_{a,b\in \Phi}$, $T=(T_{ab})_{a,b\in \Phi}$ of complex
numbers  such that:\smallskip

\item{{\bf MD1.}} $S$ is unitary and symmetric, and $T$ is diagonal
  and of finite order: i.e.\ $T^N=I$ for some $N$;

\item{{\bf MD2.}} $S_{0a}>0$ for all $a\in \Phi$;

\item{{\bf MD3.}} $S^2=(ST)^3$;

\item{{\bf MD4.}} The numbers defined by
$$N_{ab}^c=\sum_{d\in \Phi}{S_{ad}\,S_{bd}\,\overline{S_{cd}}\over S_{0d}}\eqno(
2.1)$$
are in $\Z_{\ge}$.\smallskip

The matrix $S$ is more important than $T$.
The name `modular data' is chosen because $S$ and $T$ give a representation of
the (double cover of the) modular group SL$_2(\Z)$ --- as {\bf MD3} strongly hints and  as we will see in Section IV.
Trying to remain consistent
with the terminology of RCFT, we will call (2.1) `Verlinde's
formula', the $N_{ab}^c$ `fusion coefficients', and  the $a\in\Phi$
`primaries'. The distinguished primary `$0$' is called the `identity' because
of its role in the associated fusion ring, defined below.
A possible fifth axiom will be proposed shortly, and later we will propose
refinements to {\bf MD1} and {\bf MD2}, as well as a possible 6th
axiom, but in this paper we will
limit ourselves to the consequences of {\bf MD1}--{\bf MD4}.

Modular data arises
in many places in math --- some of these will be reviewed next section.
In many of these interpretations, there is for each primary $a\in\Phi$ a function
(a `character') $\chi_a:{\Bbb H}\rightarrow\C$ which yields the matrices $S$
and $T$ as in (1.2a). Also, in many examples, to each triple $a,b,c\in\Phi$
we get a vector space ${\cal H}_{ab}^c$ (an `intertwiner space' or `multiplicity
module') with dim$({\cal H}_{ab}^c)=N_{ab}^c$, and with natural
isomorphisms between ${\cal H}_{ab}^c$, ${\cal H}_{ba}^c$, etc. In many of
these examples, we have `6j-symbols', i.e.\ for any
6-tuple $a,b,c,d,e,f\in\Phi$ we have a homomorphism $\{{a\atop d}\,
{b\atop e}\,{c\atop f}\}$ from ${\cal H}_{cd}^e\otimes
{\cal H}_{ab}^c$ to ${\cal H}_{af}^e\otimes {\cal H}_{bd}^f$ obeying several
conditions (see e.g.\
[85,33] for a general treatment).
Classically, 6j-symbols explicitly described the change between the
two natural bases of the tensor product $(L_\la\otimes L_\mu)\otimes
L_\nu\cong L_\la\otimes (L_\mu\otimes L_\nu)$
of modules of a Lie group, and our 6j-symbols are their natural extension to
e.g.\ quantum groups. Characters, intertwiner spaces, and 6j-symbols don't
play any role in this paper.

If {\bf MD2} looks unnatural, think of it in the following way. It is
easy to show (using {\bf MD1} and {\bf MD4} and Perron-Frobenius
theory [55]) that some column of $S$ is nowhere 0 and of constant phase
(i.e.\ Arg$(S_{\updownarrow b})$ is constant for some $b\in\Phi$);
{\bf MD2} tells us that it is the 0 column, and that the phase is 0
(so these entries are positive). The ratios $S_{a0}/S_{00}$ are 
sometimes called {\it q(uantum)-dimensions} (see (4.2b) below).

If {\bf MD4} looks peculiar, think of it in the following way. For
each $a\in \Phi$, define matrices $N_a$ by
$(N_a)_{bc}=N_{ab}^c$. These are usually called {\it fusion matrices}.
Then {\bf MD4} tells us these $N_a$'s are simultaneously diagonalised
by $S$, with eigenvalues $S_{ad}/S_{0d}$.

The key to modular data is equation (2.1). It should look familiar
from the character theory of finite groups: Let $G$ be any finite group, let
$K_1,\ldots,K_h$ be the conjugacy classes of $G$, and write $k_i$ for the
formal sum $\sum_{g\in K_i}g$. These $k_i$'s form a basis for the centre of
the group algebra ${\Bbb C}G$ of $G$. If we write
$$k_i\,k_j=\sum_\ell c_{ij\ell}k_\ell$$
then the structure constants $c_{ij\ell}$ are nonnegative integers, and
we obtain
$$c_{ij\ell}={\|K_i\|\,\|K_j\|\,\|K_\ell\|\over \|G\|}\sum_{\chi\in {\rm Irr}\,G}
{\chi(g_i)\,\chi(g_j)\,\overline{\chi(g_\ell)}\over \chi(e)}$$
where $g_i\in K_i$. This resembles (2.1), with $S_{ab}$ replaced
with $S_{i,\chi}=\chi(g_i)$ and the identity 0 replaced with the group
identity $e$. This formal relation between finite groups and Verlinde's
formula seems to have first been noticed in [68].

More generally, modular data is closely related to association schemes
and C-algebras (as first noted in [24], and independently in [4]),
hypergroups [90], etc. That is to
say, their axiomatic systems are similar. However, the
exploration of an axiomatic system is influenced not merely by its intrinsic
nature (i.e.\ its formal list of axioms and their logical consequences),
 but also by what are perceived by the Brethren to be its characteristic examples.
There always is a context to math.
The prototypical example of a C-algebra is the space of class
functions of a finite group  while that of modular data
corresponds to the SL$_2(\Z)$ representation associated to 
an affine Kac-Moody algebra at level $k\in\Z_\ge$ (Example 2
below). Nevertheless it can be expected that techniques and questions
from one of these areas can be profitably carried over to the
other. To give one interesting disparity, the commutative association
schemes have been  classified up to 23 vertices [59], while modular data is
known for only 3 primaries [17] (and that proof assumes
additional axioms)!  In fact we still don't have a
finiteness theorem: for a given cardinality $\|\Phi\|$, are there
only finitely many possible modular data?
But see [32] for a more sophisticated and promising approach to
modular data classification.

The matrix $T$ is fairly poorly constrained by {\bf MD1}--{\bf MD4}. Another
axiom, obeyed by Examples 1,2,3 next section (as well as any conformal
field theory [25]), can be introduced,
though it won't be adopted here:

\smallskip\item{{\bf MD5.}} For all choices $a,b,c,d\in\Phi$,
$$(T_{aa}T_{bb}T_{cc}T_{dd}T_{00}^{-1})^{N_{abcd}}=\prod_{e\in\Phi}
T_{ee}^{N_{abcd,e}}$$
where
$$N_{abcd}:=\sum_{e\in\Phi}N_{ab}^e\,N_{ce}^{d}\ ,\qquad N_{abcd,e}:=
N^e_{ab}N_{ce}^d+N_{bc}^eN_{ae}^d+N_{ac}^eN_{be}^d$$

 From {\bf MD5} can be proved that $T$ has finite order
(take $a=b=c=d$), so admitting {\bf  MD5} permits us to remove that statement
from {\bf MD1}. But it doesn't have any other interesting consequences
that this author knows  --- though perhaps it will be useful in proving
the Congruence Subgroup Property given below, or give us some finiteness result.

\smallskip Intimately related to modular data are the {\it fusion rings} $R={\cal F}(\Phi,N)$.

\medskip\noindent{{\bf Definition 2.}} A fusion ring is a commutative ring $R$ with identity 1, together with a finite
basis $\Phi$ (over $\Q$) containing 1, such that:

\smallskip \item{{\bf F1.}} The structure constants
$N_{ab}^c$ are all nonnegative;

\item{{\bf F2.}} There is a ring endomorphism $x\mapsto
x^*$ stabilising the basis $\Phi$;

\item{{\bf F3.}} $N_{ab}^1=\delta_{b,a^*}$.\smallskip 

That $x\mapsto x^*$ is an involution is clear from {\bf F3} and
commutativity of $R$. Axiom {\bf F3} and associativity of $R$ imply $N_{xy}^z=
N_{xz^*}^{y^*}$ (a.k.a.\ Frobenius reciprocity or Poincar\'e duality); hence the numbers $N_{xyz}
\eqde N_{xy}^{z^*}$ will be symmetric in $x,y,z$. Axiom {\bf F3} is equivalent to
the existence on $R$ of a linear functional `Tr' for which $\Phi$ is
orthonormal: Tr$(xy^*)=\delta_{x,y}$ $\forall x,y\in\Phi$. Then $N_{xy}^z
={\rm Tr}(xyz^*)$. 
The underlying coefficient ring was chosen to be $\Q$ here, but that choice
isn't important (except that it forces the coefficients $N_{ab}^c$ to be rational).

As an abstract algebra, $R$ is not very interesting: 
in particular, because $R$ is commutative and associative,
the fusion matrices $(N_a)_{bc}=N_{ab}^c$ pairwise commute; because of {\bf F2}, $(N_a)^t=
N_{a^*}$. Thus they are normal and can be simultaneously diagonalised. Hence
$R$ is semisimple, and will be isomorphic to a direct sum of number
fields (see Example 7 below). For example, the fusion ring for $A_1^{(1)}$
level $k$ (see (3.5c)) is isomorphic to $\oplus_d\Q[\cos(\pi {d\over k+2})]$, where $d$
runs over all divisors of $2(k+2)$ in the interval $1\le d<k+2$.
Likewise, the fusion
ring $R\otimes_\Q\C$ over $\C$ is isomorphic as a $\C$-algebra to $\C^{\|\Phi\|}$ with operations
defined component-wise.
Of course what is important for fusion rings is that they have a
preferred basis $\Phi$, unlike more familiar algebras. Incidentally, more
general fusion-like rings arise naturally in subfactors
(see Example 6 below) and nonrational logarithmic CFT (see e.g.\ [48]) so their theory
also should be developed.

We usually will be interested in the `fusion coefficients' $N_{ab}^c$
being (nonnegative) integers. Note that the identity of fusion rings
is denoted here by `1' rather than the `0' used in modular data.

Our treatment now will roughly follow that of Kawada's C-algebras as given in
[6]. The fusion matrices $N_{a}$ are linearly independent, by {\bf F3}.
Let $\underline{x_i}$, for $1\le i\le n=\|\Phi\|$,
be a basis of common eigenvectors, with eigenvalues $\ell_i(a)$. Normalise
all vectors $\underline{x_i}$ to have unit length (there remains an ambiguity of
phase which we will fix below), and let $\underline{x_1}$ be the 
Perron-Frobenius one --- since $\sum_aN_a>0$ here, we can choose $\underline{x_1}$
to be strictly positive.  Let $S$ be the matrix whose $i$th column
is $\underline{x_i}$, and $L$ the matrix $L_{ai}=\ell_i(a)$. Then $S$ is
unitary and $L$ is invertible. Note that for each $i$, the map $a\mapsto \ell_i
(a)$ defines a linear representation of $R$. That means that each column of
$L$ will be a common eigenvector of all $N_a$, with eigenvalue $\ell_i(a)$,
and hence must equal a scalar multiple of the $i$th column of $S$ (see the {\smcap Basic Fact}
in Section IV). Note that
each $L_{1i}=1$; therefore each $S_{1i}$ will be nonzero and we may uniquely
determine $S$ (up to the ordering of the columns) by demanding that each
$S_{1i}>0$. Then $L_{ai}=S_{ai}/S_{1i}$. Therefore we get (2.1).

Note though that the rows of $S$ are indexed by $\Phi$, but its columns
are indexed by the eigenvectors. Like the character table of a group,
although $S$ is a square matrix it is not at this point in our
exposition `truly square'. This simple observation will be valuable for the
paragraph after Prop.\ 1.

The involution $a\mapsto a^*$ in {\bf F2} appears in the matrix $C_l:=SS^t$: $(C_l)_{ab}
=\delta_{b,a^*}$. The matrix $C_r:=S^tS$ is also an order 2 permutation,
and 
$$\overline{S_{ai}}=S_{C_la,i}=S_{a,C_ri}\eqno(2.2)$$
For a proof of those statements, see (4.4) below.

Let $\widehat{R}$ be the set of all linear maps of $R\otimes_\Q\C$ into $\C$, equivalently
the set of all maps $\Phi\rightarrow\C$. $\widehat{R}$ has the structure of
an $(n+1)$-dimensional commutative algebra over $\C$, using the product $(fg)(a)=f(a)g(a)$.
A basis $\widehat{\Phi}$ of $\widehat{R}$ consists of the functions $a\mapsto 
{S_{ai}\over S_{a1}}$, for each $1\le i\le n$ --- denote this function $\widehat{i}$.
 The resulting structure constants are
 $$\widehat{N}_{\hat{i}\hat{j}}^{\hat{k}}
 =\sum_{a\in\Phi}{S_{ai}\,S_{aj}\,\overline{S_{ak}}\over S_{a1}}=:\widehat{N}_{ij}^k\eqno(2.3)$$
 In other words, replace $S$ in (2.1) with $S^t$.
It is easy to verify that $\widehat{R}={\cal F}(\widehat{\Phi},\widehat{N})$ obeys all axioms of a 
fusion ring (over $\C$ rather than $\Q$), except possibly that the structure constants may not be 
nonnegative. They will necessarily be real, however. We call $\widehat{R}={\cal F}(\widehat{\Phi},\widehat{N})$
the {\it dual} of $R={\cal F}(\Phi,N)$. Note that $\widehat{\hat{R}}$
can always be naturally identified with $R\otimes \C$.

We call $R={\cal F}(\Phi,N)$ {\it
self-dual} if $\widehat{R}={\cal F}(\widehat{\Phi},\widehat{N})$ is isomorphic
as a fusion ring 
to $R\otimes \C$ --- equivalently, if there is a bijection $\io:\Phi\rightarrow
\widehat{\Phi}$
 such that ${N}_{{a}{b}}^{{c}}=\widehat{N}_{\io a,\io b}^{\io c}$ (see the
definition of `fusion-isomorphism' in Section IV).

\medskip\noindent{{\smcap Proposition 1.}}\quad {\it Given any fusion ring $R={\cal F}(\Phi,N)$,
there is a unique (up to ordering of the columns) unitary matrix $S$
obeying $(2.1)$ and all $S_{1i}$ and $S_{a1}$ are positive. The fusion ring $R={\cal F}(\Phi,N)$
is self-dual iff the corresponding matrix $S$ obeys
$$S_{a,\iota' b}=S_{b,\iota a}\qquad{\rm for\ all}\ a,b\in\Phi\eqno(2.4)$$
for some bijections $\io,\io':\Phi\rightarrow\widehat{\Phi}$.}\medskip

What this tells us is that there isn't a natural algebraic
 interpretation for our condition $S=S^t$ in {\bf MD1}; the study of fusion rings insists
 that the definition of modular data be extended to the more general
setting where  `$S=S^t$' is replaced with (2.4). Fortunately,
all properties of modular data extend naturally to this new setting.
But what should $T$ look like then? A priori this isn't so
 clear. But requiring the existence of a representation of SL$_2(\Z)$ really forces
 matters. In particular note that, when $S$ is not symmetric, the matrices $S$ and $T$
 themselves cannot
 be expected to give a natural representation of any group (modular or
 otherwise) since for instance the expression $S^2$ really isn't 
sensible --- $S$ is not `truly square'. Write $P$ and
 $Q$ for the matrices $P_{a,i}=\delta_{i,\io a}$ and
 $Q_{a,i}=\delta_{i,\io' a}$, and let $n$ be the order of the
 permutation $\io^{-1}\circ \io'$. 
Then for any $k$,
 $\tilde{S}=SQ^t(PQ^t)^k$ {\it is} `truly square' and its square
 $\tilde{S}^2=C_l(P Q^t)^k$ is a permutation matrix, where $C_l$  is as in
 (2.2). We also want
 $\tilde{S}^4=I$, which requires $n=2k+1$ or $n=4k+2$. In either of those cases,
 $\tilde{T}=TP^t(QP^t)^k$ defines with $\tilde{S}$ a representation of
 SL$_2(\Z)$ provided $TS^tTS^tT=S(Q^tP)^{2k+1}$. 
(When 4 divides $n$, the best we will get
in general will be a representation of some extension of SL$_2(\Z)$.) 
But $S$ is only determined by the fusion ring up to permutation of the
columns, so we may as well replace it with $\tilde{S}$. Do likewise
 with $T$. 
So it seems that  we can and should replace {\bf MD1} with:

\smallskip\noindent{{\bf MD1$'$.}} $S$ is unitary, $S^t=SP$ where $P$
is a permutation matrix of order a power of 2, and $T$ is diagonal and
of finite order;

\smallskip\noindent and leave {\bf MD2}--{\bf MD4} intact.
That simple change seems to provide the natural generalisation of modular
 data to any self-dual fusion ring. Let $n$ be the order of $P$;
 then $n=1$ recovers modular data, $n\le 2$ yields a representation of
SL$_2(\Z)$, and $n>2$ yields a representation of a central extension of
SL$_2(\Z)$.

If we don't require an SL$_2(\Z)$ representation, then
 of course we get much more freedom.
 It is very unclear though what $T$ should look like when the fusion ring is not self-dual,
which probably indicates that  the definition
 of fusion ring should include some self-duality constraint.  
This is the attitude we adopt.

Incidentally,
the natural appearance of a self-duality constraint 
here perhaps should not be surprising in hindsight.
Drinfeld's `quantum double' construction has analogues in several contexts,
and is a way of generating algebraic structures which possess modular data
(see examples next section). It always involves combining a given (inadequate)
algebraic structure with its dual in some way. A general categorical
interpretation of quantum double is the {\it centre construction}, described
for instance in [67]; it assigns to a tensor category a braided tensor
category. It would be interesting to interpret this construction at the
more base level of fusion ring --- e.g.\ as a general way for obtaining
self-dual fusion rings from non-self-dual ones.

In Example 4 of Section III we will propose a further generalisation of modular data. In this paper however, we will restrict to the
consequences of {\bf MD1}--{\bf MD4}.

In any case,  a fusion ring with integral fusion
coefficients $N_{ab}^c$, self-dual in the strong sense that $\iota=\iota'$,
 is completely equivalent to a unitary and symmetric matrix $S$
obeying {\bf MD2}. This special case of Proposition 1 was known
to Bannai and Zuber. More generally, $\iota^{-1}\circ\iota'$ will define a
{\it fusion-automorphism} of a self-dual fusion-ring $R={\cal F}(\Phi,N)$.
Note that an unfortunate choice of matrix $S$ in
[4] led to an inaccurate conclusion there regarding fusion rings and
Verlinde's formula (2.1). In fact, Verlinde's formula will hold with a
unitary matrix $S$ obeying $S_{1i}>0$, even if we drop nonnegativity
{\bf F1}.

Proposition 1 shows that although (2.1) looks mysterious,
it is quite canonical, and that the depth of Verlinde's formula lies
in the interpretation given to $S$ and $N$ (for instance  (1.2a) and
$N_{ab}^c={\rm dim}({\cal H}_{ab}^c)$) within the given context.

The two-dimensional fusion rings ${\cal F}(\{1,2\},N)$ are classified by their value of $r=N_{22}^2$
--- there is a unique fusion ring for every $r\in\Q$, $r\ge 0$. All are
self-dual. A diagonal unitary matrix $T$ satisfying $(ST)^3=S^2$
exists, iff $0\le r\le {2\over \sqrt{3}}$. However, $T$ will in
addition be of finite order, i.e.\ $S$ and $T$ will constitute modular
data, iff $r=0$
(realised e.g.\ by the affine algebras $A_1^{(1)}$ and $E_7^{(1)}$ level 1)
or $r=1$ (realised e.g.\ by  $G_2^{(1)}$ and $F_4^{(1)}$ level 1).
Both $r=0,1$ have six possibilities for the matrix $T$ ($T$ can always be
multiplied by a third root of unity). All 12 sets of modular data with two primaries
can be realised by affine algebras (see Example 2 below). This seeming
 omnipresence of the  affine algebras is an accident of small numbers of primaries;
 even when $\|\Phi\|=3$ we find non-affine algebra modular data. The (rational)
fusion rings given here can be regarded as a deformation interpolating
between e.g.\ the $A_1^{(1)}$ and $G_2^{(1)}$ level 1 fusion rings;
similar deformations exist in higher dimensions. For example in 3-dimensions,
the $A_2^{(1)}$ level 1 fusion ring lies in a family of (rational) self-dual
fusion rings parametrised by the Pythagorean triples.

{\it Classifying modular data and fusion rings for small sets of primaries, 
or at least obtaining new explicit families beyond Examples 1-3 given next section,
is perhaps the most vital challenge in the theory.}

\bigskip\centerline{{\bf III. Examples of modular data and fusion rings}}\medskip

We can find (2.1), if not modular data in its full splendor, in a wide
variety of contexts. In this section we sketch several of these.
Historically for the subject, Example 2 has been the most important. As with the introductory section,
don't be concerned if most of these examples aren't familiar --- just
move on to Section IV.

\medskip\noindent{{\smcap Example 1:} {\it Lattices.}} See [19]
for the essentials of lattice theory.

Let $\L$ be an
even lattice --- i.e.\ $\L$ is the $\Z$-span of a basis of $\R^n$,
with the property that $x\cdot y\in \Z$ and $x\cdot x\in 2\Z$ for all 
$x,y\in \L$. Its dual $\L^*$ consists of all vectors $w$ in $\R^n$ whose
dot product $w\cdot x$ with any $x\in \L$ is an integer. So we have 
$\L\subset \L^*$.
 Let $\Phi=\L^*/\L$ be the cosets. The cardinality of $\Phi$ is finite,
 given by the determinant $|\L|$ of $\L$ (which equals the
volume-squared of any fundamental region).
 The dot products $a\cdot b$ and norms $a\cdot a$ for the classes
 $[a],[b]\in\Phi$ are well-defined (mod $1$) and (mod 2), respectively.
 Define matrices by
 $$\eqalignno{S_{[a],[b]}=&\,{1\over \sqrt{|\L|}}e^{2\pi\i a\cdot b}&(3.1a)\cr
 T_{[a],[a]}=&\,e^{\pi\i a\cdot a-n\pi\i/12}&(3.1b)\cr}$$
The simplest special case is $\L=\sqrt{N}\Z$ for any even number $N$,
 where $\L^*={1\over\sqrt{N}}\Z$ and $|\L|=N$.
Then $\Phi$ can be identified with $\{0,1,\ldots,N-1\}$, and $a\cdot b$
is given by $ab/N$, so (3.1a) becomes the finite Fourier transform.

For any such lattice $\L$, this defines modular data. Note that the
SL$_2(\Z)$-representation is essentially a Weil representation of
SL$_2(\Z/|\Lambda|\Z)$, and that it is realised in the sense of (1.2) by
characters ch$_{[a]}$ given by theta functions divided by
$\eta(\tau)^n$. The identity `0' here is
$[0]=\L$. The fusion
coefficients $N_{[a],[b]}^{[c]}$ equal the Kronecker delta 
$\delta_{[c],[a+b]}$, so the product in the fusion ring is given by addition
in $\L^*/\L$. From our point of view, this lattice example is 
too trivial to be interesting.

When $\L$ is merely integral (i.e.\ some norms $x\cdot x$ are odd), we
don't have  modular data:
$T^2$ (but not $T$) is defined by (3.1b), and we get a representation of $\langle
\left(\matrix{0&-1\cr 1&0}\right),\left(\matrix{1&2\cr 0&1}\right)\rangle$,
an index-3 subgroup of SL$_2(\Z)$. However, nothing essential is lost, so
the definition of modular data should be broadened to include at
minimum all these integral lattice examples. 
\QED

\medskip\noindent{{\smcap Example 2:} {\it Kac-Moody algebras.}}
See [63,66] for the basics of Kac-Moody algebras.

The source of some of the most interesting modular data are the
affine nontwisted Kac-Moody algebras $X_r^{(1)}$. The simplest
way to construct affine algebras is to let $X_r$ be any finite-dimensional simple
(more generally, reductive) Lie algebra. Its loop algebra is the set of all
formal series $\sum_{\ell\in \Z}t^\ell a_\ell$, where $t$ is an
indeterminant, $a_\ell\in X_r$
and all but finitely many $a_\ell$ are 0. This is a Lie algebra, using
the obvious bracket, and is infinite-dimensional. The affine algebra $X_r^{(1)}$
is simply a certain central extension of the loop algebra. (As usual,
the central extension is taken in order to get a rich supply of
representations.) 

The representation theory of $X_r^{(1)}$ is analogous to that of $X_r$.
We are interested in the integral highest weight representations. These
are partitioned into finite families parametrised by the level $k\in\Z_{\ge}$.
Write $P_+^k(X_r^{(1)})$ for the finitely many level $k$ highest weights
$\la=\la_0\L_0+\la_1\L_1+\cdots+\la_r\L_r$, $\la_i\ge 0$. For example, $P_+^k(A_r^{(1)})$
consists of the $\left({k+r\atop r}\right)$ such  $\la$, which obey 
$\la_0+\la_1+\cdots+\la_r=k$.

The $X_r^{(1)}$-character $\chi_\la(\tau)$ associated to
highest weight $\la$ is given by a graded trace, as in (1.1). Thanks
to the structure and action of the affine Weyl group on the Cartan subalgebra
of $X_r^{(1)}$, the character $\chi_\la$ is
essentially a lattice theta function, and so transforms nicely under the modular
group SL$_2(\Z)$. In fact, for fixed algebra $X_r^{(1)}$
and level $k\in\Z_{\ge}$, these $\chi_\la$ define a representation of SL$_2(\Z)$,
exactly as in (1.2) above, and the matrices $S$ and $T$
constitute modular data. The `identity' is $0=k\L_0$, and
the set of `primaries' is the highest weights $\Phi=P_+^k(X_r^{(1)})$. The matrix $T$ is related to the
values of the second Casimir of $X_r$, and $S$ to elements of finite order in the
Lie group corresponding to $X_r$:
$$\eqalignno{T_{\la\mu}=&\,\alpha\,\exp[{\pi\i\,(\la+\rho|\la+\rho)\over \k}]\,
\delta_{\la,\mu}&(3.2a)\cr S_{\mu\nu}=&\, \alpha'\,
\sum_{w\in \overline{W}}{\rm det}(w)\,\exp[-2\pi \i\,{(w(\mu+\rho)|
\nu+\rho)\over \k}]\ &(3.2b)\cr
{S_{\la\mu}\over S_{0\mu}}=&\,{\rm ch}_{\overline{\la}}(\exp[-2\pi\i\,{(\overline{\la}
\,|\,\overline{\mu}+\rho)\over\k}])&(3.2c)}$$
The numbers  $\alpha,\alpha'\in {\Bbb C}$ are normalisation constants whose precise values are unimportant
here, and are given in Thm.\ 13.8 of [63].  The inner product in (3.2)
is the usual Killing form, $\rho$ is the Weyl vector
$\sum_i\L_i$, and $\kappa=k+h^\vee$, where $h^\vee$ is the dual
Coxeter number ($=r+1$ for $A_r^{(1)}$).
The (finite) Weyl group $\overline{W}$ of $X_r$
 acts on $P_+^k$ by fixing $\L_0$.
Here, $\overline{\la}$ denotes the projection $\la_1\L_1+\cdots+\la_r\L_r$,
and `ch$_{\overline{\la}}$' is a finite-dimensional Lie group character.

The combinatorics of Lie group characters at elements of finite order,
i.e.\ the ratios (3.2c), is quite rich.
For example, in [62] they are used to prove quadratic reciprocity, while [72] 
uses them for instance in a fast algorithm for computing tensor product decompositions
in Lie groups.

The fusion coefficients $N_{\la\mu}^\nu$, defined by (2.1),
 are essentially the tensor product
multiplicities $T_{{\la}{\mu}}^{{\nu}}\eqde\,$mult$_{
\overline{\la}\otimes\overline{\mu}}(\overline{\nu})$ 
for $X_r$ (e.g.\ the Littlewood-Richardson coefficients for
$A_r$), except `folded' in a way depending on $k$. This is seen
explicitly by the Kac-Walton formula [63 p.\ 288, 88,44]:
$$N_{\la\mu}^\nu=\sum_{w\in{W}}{\rm det}(w)\,T_{{\la}{\mu}}^{
{w.\nu}}\ ,\eqno(3.3)$$
where $w.\ga\eqde w(\ga+\rho)- \rho$ and $W$ is the affine Weyl group of
$X_r^{(1)}$ (the dependence on $k$ arises through this action of $W$).
The proof of (3.3) follows quickly from (3.2c).

The fusion ring $R$ here is isomorphic to Ch$(X_\ell)/{\cal I}_k$,
where Ch$(X_\ell)$ is the character ring of $X_\ell$ (which is
isomorphic as an algebra to the polynomial algebra in $\ell$
variables), and where ${\cal I}_k$ is its ideal generated by the
characters of the `level $k+1$' weights (for $X_\ell=A_\ell$, these
consist of all $\overline{\la}=(\la_1,\ldots,\la_\ell)$ obeying
$\la_1+\cdots+ \la_\ell=k+1$).

Equation (3.3) has the flaw that, although it is manifest that the $N_{\la\mu}^\nu$
will be integral, it is not clear why they are positive.  A big open challenge
here is the discovery of a combinatorial rule, e.g.\ in the spirit of the
well-known Littlewood-Richardson rule, for the affine fusions. Three
preliminary steps in this direction are [82,84,35].

Identical numbers $N_{\la\mu}^\nu$ appear in several other contexts. For instance,
Finkelberg [36] proved that the affine fusion ring is isomorphic to the K-ring of
Kazhdan-Lusztig's category $\widetilde{\cal O}_{-k}$ of level $-k$
integrable highest weight $X_r^{(1)}$-modules, and to Gelfand-Kazhdan's
category $\widetilde{\cal O}_q$ coming from finite-dimensional modules of the
quantum group $U_q(X_r)$ specialised
to the root of unity $q=\exp[\i\pi/m\kappa]$ for appropriate choice of $m\in
\{1,2,3\}$.  Because of these isomorphisms, we know that
the $N_{\la\mu}^\nu$ do indeed lie in $\Z_\ge$, for any affine algebra.
We also know [38] that they increase with $k$, with limit $T_{\la\mu}^\nu$.

Also, they arise as dimensions of spaces of generalised theta functions [34],
as tensor product coefficients in quantum groups [44] and Hecke algebras [58] at roots
of 1 and Chevalley groups for ${\Bbb F}_p$ [56], and in quantum cohomology [91].

For an explicit example, consider the simplest affine algebra $(A_1^{(1)}$)
at level $k$.
We may take $P_+^k=\{0,1,\ldots,k\}$
(the value of $\la_1$),
and then the $S$ and $T$ matrices and fusion coefficients are given by
$$\eqalignno{S_{ab}=&\,\sqrt{{2\over k+2}}\,\sin(\pi\,{(a+1)\,
(b+1)\over k+2})&(3.5a)\cr T_{aa}=&\,\exp[{\pi\i\, (a+1)^2\over 2(k+2)}
-{\pi\i\over 4}]&(3.5b)\cr N_{ab}^c=&\,\left\{\matrix{1&{\rm if}\ c\equiv a\!+\!b
\ ({\rm mod\ 2)\ and}\ |a\!-\!b|\le c\le{\rm min}\{a\!+\!b,
2k\!-\!a\!-\!b\}\cr0&{\rm otherwise}\cr}\right. &(3.5c)}$$
The only other affine algebras for which the fusions have been
explicitly calculated are $A_2^{(1)}$ [9] and $A_3^{(1)}$ [10],
and their formulas are also surprisingly compact.

Incidentally, an analogous modular transformation matrix $S$ to (3.2b)
exists for the so-called {\it admissible representations} of $X_r^{(1)}$
at fractional level [65]. The matrix is symmetric, but has no
column of constant phase and thus naively putting it into Verlinde's
formula (2.1) will necessarily produce some negative numbers (apparently they'll always
be integers though). A legitimate fusion ring has been obtained for 
$A_1^{(1)}$ at fractional level $k={p\over q}-2$ in other ways [3]; it
factorises into the product of the $A_{1,p-2}$ fusion ring with a fusion ring at
`level' $q-1$ associated to the rank 1 supersymmetric algebra osp$(1|2)$.
Some doubt however on the relevance of these efforts has been cast by
[46]. A similar theory should
exist at least for the other $A_r^{(1)}$; initial steps for $A_2^{(1)}$
have been made in [45].
Exactly how these correspond to
modular data, or rather how modular data should be generalised to accommodate them, is
not completely understood at this time. \QED

\medskip\noindent{{\smcap Example 3:} {\it Finite groups.}} The
relevant aspects of finite group theory are given in  e.g.\ [61].

Let $G$ be any finite group. Let $\Phi$ be the set
of all pairs $(a,\chi)$, where the $a$ are representatives of the conjugacy
classes of $G$ and $\chi$ is the character of an irreducible representation of
the centraliser $C_G(a)$. (Recall that the conjugacy class of an
element $a\in G$ consists of all elements of the form $g^{-1}ag$, and
that the centraliser $C_G(a)$ is the set of all $g\in G$ commuting
with $a$.) Put [26,71]
$$\eqalignno{S_{(a,\chi), (a',\chi')}\,&={1\over \|C_G(a)\|\,\|C_G(a')\|}\sum_{g\in G(a,a')}
\overline{\chi'(g^{-1}ag)}\,\overline{\chi(ga'g^{-1})}&(3.6a)\cr 
T_{(a,\chi),(a',\chi')}\,&=\delta_{a,a'}\delta_{\chi,\chi'}{\chi(a)\over\chi(e)}&(3.6b)}$$
where $G(a,a')=\{g\in G\,|\, aga'g^{-1}=ga'g^{-1}a\}$, and $e\in G$ is the
identity. For the `identity'  0 take $(e,1)$. Then (3.6) is modular
data. See [22] for several explicit examples.

There are  group-theoretic descriptions
of the fusion coefficient $N_{(a,\chi),(b,\chi')}^{(c,\chi'')}$. 
That these fusion coefficients are nonnegative integers, follows for
instance from
Lusztig's interpretation of the corresponding fusion ring as the Grothendieck ring
of equivariant vector bundles over $G$: $\Phi$ can be identified with
the irreducible vector bundles. 

This class of modular data
played an important role in Lusztig's determination of
irreducible characters of Chevalley groups. But there is a remarkable
variety of contexts in which (3.6) appears (these are reviewed in [22]). For instance,
modular data often has a Hopf algebra interpretation: just as the
affine fusions are recovered from the quantum group $U_q(X_r)$, so are
these finite group fusions recovered from the quantum-double of $G$.

This modular data is quite interesting for nonabelian $G$, and
deserves more study. It behaves very differently than the affine data [22]. Conformal
field theory explains how very general
constructions (Goddard-Kent-Olive and orbifold) build up modular data
from combinations of affine and finite group data --- see e.g.\ [25]. 

For a given finite group $G$, there doesn't appear to be a natural
unique choice of characters ch$_{(a,\chi)}$ realising this modular
data in the sense of (1.2).

This modular data can be twisted [27] by a 3-cocycle
$\alpha\in H^3(G,\C^\times)$, which plays the same 
role here that level did in Example 2. A further major generalisation of this finite group data will be
discussed in Example 6 below, and of this cohomological twist $\alpha$
in the paragraph after Example 6. \QED

\medskip \noindent{{\smcap Example 4:} {\it RCFT, TFT.}} See e.g.\
[25] and [85], and references therein, for good surveys of
2-dimensional conformal and 3-dimensional topological field theories,
respectively. In [39] can be found a survey of fusion rings in rational
conformal field theory (RCFT).

As discussed earlier, a major source of modular data comes from RCFT (and
string theory) and, more or less the same thing, 3-dimensional topological
field theory (TFT).  

In RCFT, the elements $a\in\Phi$ are called `primary fields', and
the privileged one `0' is called the `vacuum state'.
The entries of $T$ are interpreted in RCFT to
be $T_{aa}=\exp[2\pi\i (h_a-{c\over 24})]$, where $c$ is the rank of the
VOA or the
`central charge' of the RCFT, and $h_a$ is the `conformal
weight' or $L_0$-eigenvalue of the primary field $a$.
Equation (2.1) is a special case of the so-called {\it Verlinde's formula} [87]:
$$V^{(g)}_{a^1\cdots a^t}=\sum_{b\in \Phi} (S_{0b})^{2(1-g)}{S_{a^1b}\over
S_{0b}}\cdots {S_{a^tb}\over S_{0b}}\eqno(3.7)$$
It arose first in RCFT as an extremely useful expression for
the dimensions of the space of conformal blocks on a genus $g$ surface
with $t$ punctures, labelled with primaries $a^i\in\Phi$ --- the fusions 
$N_{ab}^c$ correspond to a sphere with 3 punctures.
All the $V^{(g)}$'s are nonnegative integers iff all the $N_{ab}^c$'s
are.  In RCFT, our unused axiom  {\bf MD5}  
is derived by applying Dehn twists to a sphere with 4 punctures to obtain an
$N_{abcd}\times N_{abcd}$ matrix equation on the corresponding space of conformal blocks;
{\bf MD5} is the determinant of that equation [86].

Example 1 corresponds to the string theory of $n$ free bosons compactified on the
torus $\R^n/\L$. Example 2 corresponds to Wess-Zumino-Witten RCFT [57]
where a closed string lives on a Lie group manifold. Example 3 corresponds to
the untwisted sector in an orbifold of a holomorphic RCFT (a holomorphic theory has trivial
modular data --- e.g.\ a lattice theory when the lattice $\L=\L^*$ is self-dual) by $G$ [26]. The RCFT
interpretation of fractional level affine algebra modular data isn't 
understood yet, despite considerable effort (see e.g.\ [46]). 

Actually, in a {\it nonunitary} RCFT, the matrices $S$ and
$T$ defined by (1.2a) will obey {\bf MD1}, {\bf MD3}, and {\bf MD4},
but not {\bf MD2}. For example, the `$c=c(7,2)=-{68\over 7}$
nonunitary minimal model' has $S$ and $T$, defined by (1.2a), given by
$$\eqalignno{T=&\,{\rm diag}\{\exp[17\pi\i/21],\exp[5\pi\i/21],\exp[-\pi\i/21]\}&\cr
S=&\,{2\over \sqrt{7}}\left(\matrix{\sin(2\pi/7)&-\sin(3\pi/7)&\sin(\pi/7)\cr
-\sin(3\pi/7)&-\sin(\pi/7)&\sin(2\pi/7)\cr \sin(\pi/7)&\sin(2\pi/7)&
\sin(3\pi/7)}\right)&(3.8a)}$$
This is not modular data, since the first column is not strictly
positive. However the 3rd column is. The nonunitary RCFTs tell us to 
 replace {\bf MD2} with\smallskip 

\noindent{\bf MD2$'$.} For all $a\in\Phi$, $S_{0,a}$ is a nonzero real
number. Moreover there is some $0'\in\Phi$ such that $S_{0',a}>0$ for
all $a\in\Phi$.\smallskip
 
Incidentally, an $S$ matrix which the algorithm of Section II
would  associate to that $c=-{68\over 7}$ minimal model is
$$S={2\over \sqrt{7}}\left(\matrix{\sin(\pi/7)&\sin(2\pi/7)&\sin(3\pi/7)\cr
\sin(2\pi/7)&-\sin(3\pi/7)&\sin(\pi/7)\cr\sin(3\pi/7)&\sin(\pi/7)&
-\sin(2\pi/7)}\right)\eqno(3.8b)$$
We can tell by looking at (3.8b) that it can't directly be given
the familiar interpretation (1.2a). The reason is that any such matrix
$S$ must have a strictly positive eigenvector with eigenvalue 1:
namely the eigenvector with $a$th component ${\rm ch}_a(\i)$
($\tau=\i$ corresponds to $q=e^{-2\pi}>0$ and is fixed by
$\tau\mapsto -1/\tau$; moreover the characters of VOAs converge at any
$\tau\in{\Bbb H}$ [92]). Unlike the $S$ in (3.8a), the $S$ of (3.8b)
has no such eigenvector. Thus we may find it convenient (especially
in classification attempts) to introduce a new axiom:

\smallskip\noindent{{\bf MD6.}} $S$ has a strictly positive
eigenvector $\underline{x}>0$ with eigenvalue 1.\smallskip

Note that with the choice $T={\rm
diag}\{\exp[\pi\i/21],\exp[-17\pi\i/21],\exp[-5\pi\i/21]\}$, (3.8b)
obeys {\bf MD1}-{\bf MD4}. Remarkably, all nonunitary RCFT known to
this author behave similarly. In fact, the following refinement of
{\bf MD2}$'$ appears to be true: The primary $0'$ in {\bf MD2}$'$
equals $J\sigma 0$ for some simple-current $J$ and some Galois
automorphism $\sigma$. (The term `simple-current' and the Galois
action on $\Phi$ will be defined in Section IV.)  For example, in
(3.8a) the simple-current is the identity and the Galois automorphism
corresponds to $5\in(\Z/42\Z)^\times$. Whenever this refinement holds
(which may be always), one consequence will be that the fusion ring
is realised by modular data satisfying {\bf MD1}-{\bf MD4}.

Knot and link invariants in $S^3$ (equivalently,
$\R^3$) can be obtained from an $R$ matrix and
braid group representations --- e.g.\ we have this with any quasitriangular
Hopf algebra. The much richer structure of 
{\it topological field theory} (or, 
in category theoretic language, a {\it modular category} [85]) gives us
link invariants in any closed 3-manifold, and with it modular data. 
In particular, the $S$ entries correspond to the invariants of the
Hopf link in $S^3$, $T$ to the eigenvalues of the twist operation
(Reidemeister 1, which won't act trivially here --- strictly speaking,
we have knotted ribbons, not strings), and the fusion coefficients to the invariants of 3
parallel circles $S^1\times \{p_1,p_2,p_3\}$ in the manifold
$S^1\times S^2$. Link invariants are obtained for arbitrary closed
3-manifolds by performing 
Dehn surgery,
transforming the manifold  into $S^3$; the condition that the resulting
invariants be well-defined, independent of the specific Dehn moves which get us
to $S^3$, is essentially the statement that $S$ and
$T$ form a representation of SL$_2(\Z)$. This is all discussed very
clearly in [85]. For instance, we get $S^3$ knot invariants from the
quantum group $U_q(X_r)$ with generic parameter,
but to get modular data requires specialising $q$ to a root of unity.

For extensions of this picture to representations of higher genus mapping
class groups, see e.g.\ [43] and references therein, but there
is much more work to do here. \QED

\medskip \noindent{{\smcap Example 5:} {\it VOAs.}} 
See e.g.\ [37,64] for the basic facts about VOAs; the review article
[46] illustrates how VOAs naturally arise in CFT.

Another very general
source of modular data comes from vertex operator algebras (VOAs), a
rich algebraic structure first introduced by Borcherds [13]. In particular,
let ${\cal V}$ be any `rational' VOA (see e.g.\ [92] ---
actually, VOA theory is still sufficiently undeveloped that we don't yet have
a generally accepted definition of  rational VOA). Then
${\cal V}$ will have finitely many irreducible modules $M$, one of which can be identified
 with ${\cal V}$. Zhu [92] showed that their characters ch$_M(\tau)$ transform
nicely under SL$_2(\Z)$ (as in (1.2a)).
Defining $S$ and $T$ in that way,
and calling $\Phi$ the set of irreducible $M$ and the `identity'
$0={\cal V}$, we get some of the properties of modular data.

A natural conjecture is that a large class (all?) of rational VOAs
possess (some generalisation of) modular data. We know what the fusion
coefficients mean (dimension of the space of intertwiners between the
appropriate VOA modules), and what $S$ and $T$ should mean. We know that
$T$ is diagonal and of finite order, and that $S^2=(ST)^3$ is an order-2 permutation
matrix. A Holy Grail
of VOA theory is to prove (a generalisation of) Verlinde's formula for a large class of rational
VOAs. A problem is that we still don't know when (2.1) here is even defined (i.e.\ whether
all $S_{0,M}\ne 0$). However,
suppose ${\cal V}$ has the additional (natural) property that any
irreducible module $M\ne {\cal V}$ has
positive conformal weight $h_M$ ($h_M-c/24$ is the smallest power of $q$ in the
Fourier expansion of the (normalised) character ch$_M(\tau)=q^{-c/24}\sum_{n=0}^\infty a_n^Mq^{n+h_{M}}$).
This holds for instance in all VOAs associated to unitary RCFTs.
Then consider the behaviour of ch$_M(\tau)$ for $\tau\rightarrow 0$ along the
positive imaginary axis: since each Fourier
coefficient $a_n^M$ is a nonnegative number, ch$_M(\tau)$ will go to $+\infty$. But
this is equivalent to considering the limit of $\sum_{N}S_{MN}\,{\rm ch}_N(\tau)$
as $\tau\rightarrow\i\infty$ along the positive imaginary axis.
By hypothesis, this latter limit is dominated by $S_{M0}\,a_0^{0}q^{-c/24}$,
at least when $S_{M0}\ne 0$. So what we find is that, under this hypothesis,
 the 0-column of $S$ consists of nonnegative real numbers (and also that the 
rank $c$ is positive).

In this context, Example 1 corresponds to the VOA associated to the lattice
$\L$ [28]. Example 2
is recovered by [38], who find a VOA structure on the highest weight
$X_r^{(1)}$-module
$L(k\L_0)$; the other level $k$ $X_r^{(1)}$-modules $M=L(\la)$ all
have the structure of VOA modules of ${\cal V}:=L(k\L_0)$. Example 3 arises
for example in the orbifold of  a self-dual
lattice VOA by a subgroup $G$ of the automorphism group of $\L$ (see e.g.\
[31]).
An interpretation of fractional level affine algebra data could be possible
along the lines of [30], who did it for $A_1^{(1)}$ (but see [46]). \QED

\medskip\noindent{{\smcap Example 6:} {\it Subfactors.}} See e.g.\
[33,12] for good reviews of the subfactor $\leftrightarrow$ CFT relation.

 The final general
source of modular data which we will discuss comes from subfactor theory.
To start with, let $N\subset M$ be an inclusion of II$_1$ factors with finite Jones index $[M:N]$.
Even though $M$ and $N$ will often be isomorphic as factors, Jones showed that
there is rich combinatorics surrounding how $N$ is embedded in $M$. 
Write $M_{-1}=N\subset M=M_0\subset M_1\subset\cdots$ for the tower arising
from the `basic construction'. Let $\Phi_M$ denote the set of equivalence
classes of irreducible $M-M$ submodules of $\oplus_{n\ge 0}\, {}_ML^2(M_n)_M$,
and $\Phi_N$ that for the irreducible $N-N$ submodules of $\oplus_{n\ge -1}
\,{}_NL^2(M_n)_N$. Write ${\cal H}_{AB}^C$ for the intertwiner space
{Hom}$_{M-M}(C,A\otimes_MB)$. For any $A,B\in\Phi_M$, the Connes' relative tensor product
$A\otimes_M B$ can be decomposed into a direct sum $\sum_{C\in\Phi_M}
N_{AB}^CC$, where $N_{AB}^C={\rm dim}\,{\cal H}_{AB}^C\in\Z_{\ge}$
are the multiplicities. The identity is the bimodule 
${}_ML^2(M)_M$. Assume in addition
that $\Phi_M$ is finite (i.e.\ that $N\subset M$ has `finite depth').
Then all axioms of a
fusion ring will be obeyed, except possibly commutativity: unfortunately in general
$A\otimes_MB\not\cong B\otimes_MA$. 

We are interested in $M$ and $N$ being hyperfinite. An intricate subfactor 
invariant called a {\it paragroup} (see e.g.\ [75,33]) can be formulated
in terms of 6j-symbols and fusion rings [33], and resembles exactly
solvable lattice models in statistical mechanics. One way to get
modular data is by passing from $N\subset M$ to the asymptotic inclusion $\langle M,M'\cap M_\infty
\rangle\subset M_\infty$;
its paragroup will essentially be an RCFT. Asymptotic inclusion plays the
role of quantum-double here, and corresponds physically to taking the
continuum limit of the lattice model, yielding the CFT from the underlying
statistical mechanical model. More recently [76], Ocneanu has significantly
refined this construction, generalising 6j-symbols to what are called
Ocneanu cells, and extending the context to subparagroups. His new cells
have  been interpreted by [77] in terms of Moore-Seiberg-Lewellen data [73,70].

A very similar but simpler theory has been developed for type III factors. Bimodules
now are equivalent to `sectors', i.e.\ equivalence classes of endomorphisms
$\la:N\rightarrow N$ (the corresponding subfactor is $\la(N)\subset N$).
This use of endomorphisms is the key difference (and simplification) between
the type II and type III fusion theories. Given $\la,\mu\in {\rm End}(N)$, we define $\langle
\la,\mu\rangle$ to be the dimension of the vector space of intertwiners,
i.e.\ all $t\in N$ such that $t\la(n)=\mu(n)t$ $\forall n\in N$. The
endomorphism $\la\in{\rm
End}(N)$ is irreducible if $\langle \la,\la\rangle=1$. Let $\Phi={}_N\chi_N$ be a finite
set of irreducible sectors. The fusion product is given by composition
$\la\circ\mu$; addition can also be defined, and the fusion coefficient
$N_{\la\mu}^\nu$ will then be the dimension $\langle\la\mu,\nu\rangle$.
The `identity' 0 is the identity $id_N$. Restricting to a finite set $\Phi$
of irreducible sectors, closed under fusion in the obvious way, the result is similar to a fusion
ring, except again it is not necessarily commutative (after all, why
should the compositions $\la\circ\mu$ and $\mu\circ\la$ be related). The missing ingredients
are nondegenerate braidings $\epsilon^{\pm}(\la,\mu)\in{\rm Hom}(\la\mu,\mu\la)$,
which say roughly that $\la$ and $\mu$ nearly commute (the $\epsilon^\pm$
must also obey some compatibility conditions, e.g.\ the Yang-Baxter equations).
Once we have  a
nondegenerate braiding, Rehren [78] proved that we will automatically have
modular data.

We will return to subfactors in Section V. It is probably too
optimistic to hope to see in the subfactor picture to what the
characters (1.1) correspond --- different VOAs or RCFTs can correspond
it seems to equivalent subfactors. To give a simple example, the VOA
associated to any self-dual lattice will correspond to the trivial
subfactor $N=M$, where $M$ is the unique hyperfinite II$_1$
factor. With this in mind, it would be interesting 
to find an $S$ matrix arising here which violates axiom {\bf MD6}
given earlier, or the Congruence Subgroup Property of Section IV.

Jones and Wassermann have explicitly constructed the affine algebra subfactors
(both type II and III) of Example 2, at least for $A_r^{(1)}$, and
Wassermann and students Loke and Toledano Laredo later showed that they recover
the affine algebra fusions (see e.g.\ [89] for a review). Also,
to any subgroup-group pair $H<G$, we can obtain a subfactor $R\crprod H
\subset R\crprod G$ of crossed products, where $R$ is the type II$_1$ hyperfinite factor,
 and thus a (not necessarily commutative) fusion-like ring [69]. This
subfactor $R\crprod H\subset R\crprod G$ can be thought of as giving a
grouplike interpretation to $G/H$ even when $H$ is not normal.
Sometimes it will have a braiding --- e.g.\
the diagonal embedding $G<G\times G$ recovers the finite group data of Example 3.
What is intriguing is that some other pairs $H<G$ probably also have a braiding,
generalising Example 3. There is a general suspicion, due originally perhaps
to Moore and Seiberg [73] and in the spirit of Tannaka-Krein duality,
that RCFTs can always be constructed in standard ways (Goddard-Kent-Olive
cosets and finite group orbifolds) from lattice and affine algebra models.
These crossed product subfactors could conceivably provide reams of counterexamples,
suggesting that the orbifold construction can be considerably generalised.
\QED\medskip

A uniform construction of the affine algebra and finite group modular
data is provided in [27] where a 3-dimensional TFT is associated to
any topological group $G$ ($G$ will be a compact Lie group in the
affine case; $G$ is given discrete topology in the finite case). There
we see that the level $k$ and twist $\alpha$ both play the same role,
and are given by a cocycle in $H^3(G,\C^\times)$. Crane-Yetter [23] are
developing a theory of cohomological `deformations' of modular data (more precisely,
of modular categories). In [23] they discuss the infinitesimal
deformations of tensor categories,
where the objects are untouched but the arrows are deformed, 
though their ultimate interest would be in global deformations and in
particular in specialising to the especially interesting ones --- much as we
deform the enveloping algebra $U({\frak g})$ to get the quantum group
$U_q({\frak g})$ and then specialise to roots of unity to get e.g.\ 
modular data. Their work is still in preliminary stages and it
probably needs to be generalised further (e.g.\ they don't seem to recover
the level of affine algebras), but it looks very promising.
Ultimately it can be hoped that some discrete $H^3$ group will be identified
which parametrises the
different quantum doubles of a given symmetric tensor category.

Incidentally, the fact that $H^3(G,\C^\times)$ is a group strongly suggests
that it should be meaningful to compare the modular data for different
cocycles --- e.g.\ to fix the affine algebra and vary $k$. This idea
still hasn't been seriously exploited (but e.g.\ see `threshold level' in [9,10]).

There are many examples of `pseudo-modular data'.
These are interesting for probing the question of just what should be the definition of fusion
ring or modular data. Here is an intriguing example, inspired by (4.4) below.

\medskip\noindent{{\smcap Example 7} [52]: {\it Number fields.}} A basic
introduction to algebraic number theory is provided by e.g.\ [18].

Choose any finite normal extension ${\Bbb L}$ of $\Q$, and find any totally
positive $\alpha\in{\Bbb L}$ with Tr$(|\alpha|^2)=1$ (total positivity
will turn out to be necessary for {\bf F1}). Now find any $\Q$-basis
$x_1=1,x_2,\ldots,x_n$ of a subfield ${\Bbb K}$ of ${\Bbb L}$, where
$n=$deg$({\Bbb K})$, the $x_i$ being orthonormal with respect to the trace $\langle x,y\rangle_\alpha
:={\rm Tr}(|\alpha|^2x\bar{y})$ (orthonormality will guarantee {\bf
F3} to be satisfied). Let $G$ denote the set of $n$ distinct embeddings
${\Bbb K}\rightarrow \C$. Our construction requires complex conjugation to commute with all
embeddings. Under these conditions  $|\alpha|^{-2}=\sum_i |x_i|^2$.
Then we get a fusion-like ring with primaries $\Phi=\{x_1,\ldots,x_n\}$,  `$*$'
given by complex conjugation, and structure constants $N_{ij}^k={\rm Tr}
(|\alpha|^2x_ix_j\overline{x_k})\in\Q$ given by ordinary multiplication and
addition: $x_ix_j=\sum_kN_{ij}^k\,x_k$. Call
the resulting fusion-like ring ${\Bbb K}(\Phi)$.

It is easy to see that all the properties of a fusion ring are
satisfied, except possibly $N_{ij}^k\in\Q_{\ge}$.
The fusion coefficients $N_{ij}^{k}$ will
be integers  iff the $\Z$-span of the $x_i$ form an 
`{order}' of ${\Bbb K}$. We also find that the matrix
$S_{ig}=g(\alpha x_i)$,
for $g\in G$ (lift each $g$ arbitrarily to ${\Bbb L}$), diagonalises
 these fusion matrices $N_{x_i}$. This matrix $S$ is unitary, but (unless
${\Bbb K}$ is an abelian extension of $\Q$) the dual fusions $\widehat{N}$ in (2.3)
won't be rational.

Positivity {\bf F1} requires
one of the columns of $S$ to be positive; permuting with $g$, we may
require all basis elements $x_i>0$. Hence ${\Bbb K}(\Phi)$ will have
 a chance of being a
fusion ring only when ${\Bbb K}$ is `{totally real}'.

Incidentally this example is more general than it looks: it is easy to
show

\medskip\noindent{{\smcap Proposition 2.}} {\it Let $R$ be a fusion
ring which is isomorphic as a $\Q$-algebra to a field ${\Bbb K}$. Then
$R$ is isomorphic as a fusion ring to some ${\Bbb K}(\Phi)$.}\medskip

More generally, recall that an
arbitrary fusion ring (over $\Q$) is isomorphic as an algebra to a
direct sum of number fields. So an approach to studying fusion rings could
be to study how they are built up from number fields. It would be very
interesting to classify all (not necessarily self-dual) fusion rings which are isomorphic as
a $\Q$-algebra to a field.  For example, take ${\Bbb K}=\Q[\sqrt{N}]$, where
$N$ is not a perfect square, and where also any prime divisor 
$p\equiv -1$ (mod $4$) of $N$ occurs with even multiplicity. Then
we can find positive integers $a,b$ such that $N=a^2+b^2$. Take $\Phi=\{1,{b\over a}+
{1\over a}\sqrt{N}\}$, then ${\Bbb K}(\Phi)$ is a fusion ring with
$N_{22}^2={2b\over a}$.
Note that this construction exhausts all 2-dimensional rational fusion rings, except
when $\sqrt{(N_{22}^2)^2+4}$ is rational,
which corresponds to the $\Q$-algebra $\Q\oplus\Q$ (e.g.\ the fusion ring of
$A_1^{(1)}$ level 1).
For $N=5$ and $a=2$, we recover the  fusion ring
of $F_4^{(1)}$ or $G_2^{(1)}$ level 1.   \QED

\bigskip{\centerline{\bf IV. Modular data: basic theory}}\medskip

In this section we sketch the basic theory of modular data.

It is important to reinterpret (2.1) in matrix form. For each $a\in \Phi$, 
define the {\it fusion matrix} $N_a$ by
$$(N_a)_{b,c}=N_{ab}^c\ .$$
Then (2.1) says that the $N_a$ are simultaneously diagonalised by $S$.
More precisely, the $b$th column $S_{\updownarrow,b}$ of $S$ is an eigenvector of each
$N_{a}$, with eigenvalue ${S_{ab}\over S_{0b}}$. Unitarity of $S$
tells us:
${S_{ab}\over S_{0b}}={S_{ac}\over S_{0c}}$ holds  for all $a\in \Phi$,
iff $b=c$.
In other words:

\smallskip \noindent{{\smcap Basic Fact.}}
{\it All simultaneous eigenspaces are of dimension 1, and are
spanned by each column $S_{\updownarrow,b}$ .}\smallskip

Take the complex 
conjugate of (2.1): we find that $\overline{S}$ also simultaneously diagonalises the
fusion matrices $N_a$. Hence there is some permutation of $\Phi$, which we will
denote by $C$ and call {\it conjugation}, and some complex numbers $\alpha_b$, such that 
$$\overline{S_{ab}}=\alpha_b \,S_{a,Cb}\ .$$
Unitarity forces each $|\alpha_b|=1$.
Looking at $a=0$ and applying {\bf MD2}, we see that the $\alpha_b$ must be
positive.  Hence
$$\overline{S_{ab}}=S_{a,Cb}=S_{Ca,b}\eqno(4.1)$$
and so $C=S^2$. The conjugation $C$ is trivial iff $S$ is real.
Note also that $C$, like complex conjugation, is an involution, and that
$C_{00}=1$.
Some easy formulae are $N_0=I$, $N_{ab}^0=C_{ab}$, and
$N_{Ca,Cb}^{Cc}=N_{ab}^c$. Because $C=S^2=(ST)^3$, $C$ commutes with both
$S$ and $T$: $S_{Ca,Cb}=S_{a,b}$ and $T_{Ca,Cb}=T_{a,b}$.

 For example, in Example 1,
$C[a]=[-a]$, while for $A_1^{(1)}$ the matrix $S$ is real and so
$C=I$. More generally, for the affine algebra $X_r^{(1)}$ the conjugation
 $C$ corresponds to a symmetry of
the Dynkin diagram of $X_r$. For finite groups (Example 3), $C$ takes
$(a,\chi)$ to $(a^{-1},\overline{\chi})$. In RCFT, $C$ is called {\it
charge-conjugation}; it's a symmetry in quantum field theory which interchanges
particles with their antiparticles (and so reverses charge, hence the name).

Because $C$ is an involution, we know  that the assignment (1.2b)
defines a finite-dimensional representation of SL$_{2}(\Z)$, for any choice of
modular data --- hence
the name.  A surprising fact is that this representation usually (always?)
seems to factor through a congruence subgroup. We'll return to this at the
end of this section.

Perron-Frobenius theory, i.e.\ the spectral theory of nonnegative matrices
(see e.g.\ [55]), has some immediate consequences.
By {\bf MD2} and our {\smcap Basic Fact}, the Perron-Frobenius eigenvalue
of $N_a$ is ${S_{a0}\over S_{00}}$; hence we obtain the important inequality
$$S_{a0}S_{0b}\ge |S_{ab}|\,S_{00}\ .\eqno(4.2a)$$
Unitarity of $S$ applied to (4.2a) forces
$${\rm min}_{a\in\Phi} S_{a0}=S_{00}\ .\eqno(4.2b)$$
In other words the {\it q-dimensions}, defined to be the ratios ${S_{a0}\over S_{00}}$, 
are bounded below by 1. The name `q-dimension' comes from quantum
groups (and also affine algebras (3.2c)), where one finds
 a q-deformed Weyl dimension formula. In RCFT, ${S_{a0}\over S_{00}}={\rm lim}_{\tau
\rightarrow 0^+\i}{{\rm ch}_a(\tau)\over {\rm ch}_0(\tau)}$.
In the subfactor picture (Example 6), the Jones index is the square of the
q-dimension.

Cauchy-Schwarz and unitarity,
together with (4.2a), gives us the curious inequality
$$\sum_{e\in \Phi} N_{ac}^e N_{bd}^e\le {S_{a0}\over S_{00}}
{S_{b0}\over S_{00}}\eqno(4.2c)$$
for all $a,b,c,d\in\Phi$. So for instance $N_{ab}^c\le{\rm min}\{{S_{a0}\over
S_{00}},{S_{b0}\over S_{00}},{S_{c0}\over S_{00}}\}$. Equality holds in (4.2c) only if $S_{a0}=S_{b0}
=S_{00}$ (i.e.\ only if $a$ and $b$ are {\it units} --- see below).
Other inequalities are possible, though perhaps not useful: e.g.\ 
 H\"older gives us for all $a\in\Phi$ and $k,m=1,2,3,\ldots$ the following bounds on traces
 of powers of fusion matrices:
$$({\rm Tr}\,(N_a^k))^m\le \|\Phi\|^{m-1}\,{\rm Tr}\,(N_a^{km})\eqno(4.2d)$$

The inequality (4.2b) suggests that we look at those primaries $a\in\Phi$ obeying
the equality $S_{a0}=S_{00}$. Such primaries are called {\it simple-currents}
in RCFT parlance (see e.g.\ [81,25] and references therein), but the much
more obvious mathematical name is {\it units}. To any unit $j\in \Phi$, there is a phase
$\varphi_j:\Phi\rightarrow\C$ and a permutation $J$ of $\Phi$ such that
$j=J0$ and
$$\eqalignno{S_{Ja,b}&\,=\varphi_j(b)\,S_{a,b}&(4.3a)\cr
T_{Ja,Ja}\overline{T_{aa}}&\,=\overline{\varphi_j(a)}\,T_{jj}\overline{T_{00}}
&(4.3b)\cr (T_{jj}\overline{T_{00}})^2&\,=\overline{\varphi_j(j)}&(4.3c)}$$
Moreover, if $J$ is order $n$, then $\varphi_j(a)$ is an $n$th root of unity and
$T_{Ja,Ja}\overline{T_{aa}}$ is a $2n$th root of 1; when  $n$ is odd,
the latter will in fact be an $n$th root of 1.
To reflect the physics heritage, the permutation $J$
corresponding to a unit $j\in\Phi$ will be called  a simple-current.
The set of all simple-currents or units forms an abelian group (using composition of
the permutations), called the {\it centre}
of the modular data. Note that $CJ=J^{-1}C$, and
$N_{Ja,J'b}^{JJ'c}=N_{ab}^{c}$ for any simple-currents $J,J'$.

For instance, for a lattice $\L$, all $[a]\in\Phi$ are units,
corresponding to permutation $J_{[a]}([b])=[a+b]$  and phase $\varphi_{[a]}
([b])=e^{2\pi\i a\cdot
b}$. For the affine algebra $A_1^{(1)}$ at level $k$ (recall (3.5)), there is precisely one nontrivial unit,
namely $j=k$, corresponding to $J(a)=k-a$
and $\varphi_j(a)=(-1)^a$. More generally, to any affine algebra
(except for $E_8^{(1)}$ at $k=2$), the units correspond to
symmetries of the extended Dynkin diagram. For $A_1^{(1)}$ this symmetry
interchanges the 0th and 1st nodes, i.e.\ $J(\la_0\L_0+\la_1\L_1)=
\la_1\L_0+\la_0\L_1$ (recall $a=\la_1$); for $A_r^{(1)}$ the centre is
$\Z/(r+1)\Z$. In the finite group modular
data, the units are the pairs $(z,\psi)$ where $z$ lies in the centre
$Z(G)$ of  $G$, and $\psi$ is a dimension-1 character of $G$. It
corresponds to simple-current $J_{(z,\psi)}(a,\chi)=(za,\psi\chi)$ and
phase
$\varphi_{(z,\psi)}(a,\chi)=\overline{\psi(a)}\,\overline{\chi(z)}/\chi(e)$.
The centre of this modular data will thus be isomorphic to the direct
product $Z(G)\times (G/G')$, where $G'=\langle ghg^{-1}h^{-1}\rangle$
is the commutator subgroup of $G$.

To see (4.3a), note first that (4.2a) tells us $S_{0b}\ge |S_{jb}|$ for any
unit $j$, and any $b\in\Phi$. However, unitarity then forces
$S_{0b}=|S_{jb}|$, i.e.\ (4.3a) holds for $a=0$ (with $J0$ defined to be $j$),
and some numbers $\varphi_j(b)$ with modulus 1. Putting this into (2.1), we get
$N_jN_{Cj}=I$,
the identity matrix. But the only nonnegative integer matrices whose inverses
are also nonnegative integer matrices, are the permutation matrices. This
defines the permutation $J$ of $\Phi$. Equation (4.3a) now follows from
Cauchy-Schwartz  applied to
$$1=N_{j,a}^{Ja}=\sum_{d\in\Phi}\varphi_j(d)\,S_{ad}\,\overline{S_{Ja,d}}$$
The reason $J\circ J'=J'\circ J$ is because the fusion matrices
commute: $N_{J\circ J'}=N_JN_{J'}=N_{J'}N_J=N_{J'\circ J}$.

To see (4.3b), first write $(ST)^3=C$ as $STS=\overline{T}S\overline{T}$,
then use that and (4.3a) to show $(\overline{T}S\overline{T})_{Ja,0}
=(\overline{T}S\overline{T})_{a,J0}$. To see (4.3c), use (4.3b) with $a=J^{-1}0$,
together with the fact that $C$ commutes with $T$. Note that $\varphi_j(j')=S_{j,j'}
/S_{00}=\varphi_{j'}(j)$ and $\varphi_{JJ'0}(a)=\varphi_j(a)\,\varphi_{j'}(a)$,
so $\varphi_{J^k0}(a)=(\varphi_{j}(a))^k$; from all these and (4.3b) we get
that 
$$1=T_{J^n0,J^n0}\,\overline{T_{00}}=\overline{\varphi_j(j)}\,{}^{n(n-1)/2}
(T_{jj}\overline{T_{00}})^n$$

Equations (4.2a) and (4.3b) also follow from the curious equation
$$\overline{S_{ab}}\,T_{aa}\,T_{bb}\,\overline{T_{00}}=\sum_{c\in\Phi}
N_{ab}^c\,T_{cc}\,S_{c0}$$
which is derived from (2.1) and $STS=\overline{T}S\overline{T}$.

Simple-currents and units play an important role in the theory of modular data and fusion
rings. One place they appear
is gradings.
By a {\it grading} on $\Phi$ we mean a map $\varphi:\Phi\rightarrow
\C^\times$ with the property that if $N_{ab}^c\ne 0$
then $\vp(c)=\vp(a)\,\vp(b)$. 
The phase $\varphi_j$ coming from a unit is clearly a grading;
a little more work [52] shows that any grading $\varphi$ of $\Phi$ corresponds to
a unit $j$ in this way.  The multiplicative group of gradings,
and the group of simple-currents (the centre), are naturally isomorphic.

Next, we will generalise the conjugation
symmetry argument, to other  Galois automorphisms. In particular,
write ${\Bbb Q}[S]$ for the field generated over ${\Bbb Q}$ by all entries
$S_{ab}$. Then for
 any Galois automorphism $\si\in {\rm Gal}({\Bbb Q}[S]/{\Bbb Q})$,
$$\si(S_{ab})=\epsilon_\si(a)\,S_{\si a,b}=\epsilon_\si(b)\,S_{a,\si b}\eqno(4.4)$$
for some permutation $c\mapsto \si c$ of $\Phi$, and some signs 
$\epsilon_\si:\Phi\rightarrow\{\pm 1\}$. Moreover,  the
complex numbers $S_{ab}$ will necessarily lie in the cyclotomic
extension $\Q[\xi_n]$ of $\Q$, for some root of unity $\xi_n:=\exp[2\pi\i/n]$.

For a field extension ${\Bbb K}$
of $\Q$, Gal$({\Bbb K}/\Q)$ denotes the automorphisms $\sigma$ of ${\Bbb K}$
fixing all rationals. 
Recall that each automorphism $\si\in{\rm Gal}(\Q[\xi_n]/\Q)$ corresponds to an integer 
$1\le \ell\le n$ coprime to $n$, acting by $\si(\xi_n)=\xi_n^\ell$. Note that equation (4.4)
tells us the power $\si^{2\|\Phi\|!}$ will act trivially on each entry $S_{ab}$. In
other words, the degree of the field extension $[\Q[S]:\Q]$ is bounded
by (in fact divides) $2\|\Phi\|!$. This is perhaps the closest we have
to a finiteness result for modular data (see however [7] which obtains a
bound for $n$ in terms of $\|\Phi\|$, for the modular data arising in
RCFT).

In other
incarnations, this Galois action appears in the $\chi(g)\mapsto \chi(g^\ell)$ symmetry
of the character table of a finite group, and of the action of SL$_2(\Z/N\Z)$
on level $N$ modular functions.
Equation (4.4) was first shown in [20] and a related
symmetry for commutative association schemes was found in [74]. 
The analogue of cyclotomy isn't known for  association schemes. The reason is
the additional `self-duality' property of the fusion ring, i.e.\ the fact that
$S=S^t$ or more generally (2.4).

Recall from Section II that a fusion ring $R={\cal F}(\Phi,N)$ is
isomorphic to a direct sum of number fields.
The Galois orbits determine these fields. In particular, for any
Galois orbit $[d]$
in $\Phi$, let ${\Bbb K}_{[d]}$ denote the field generated by all numbers
of the form ${S_{ab}\over S_{0b}}$ for $a\in\Phi$ and $b\in[d]$. Then
$R$ is isomorphic as a $\Q$-algebra to the direct sum
$\oplus_{[d]}{\Bbb K}_{[d]}$. We gave the $A_1^{(1)}$ level $k$
example in Section II.

The Galois action for the lattice modular data is simple: the Galois automorphisms
$\si=\si_\ell$ correspond to integers $\ell$ coprime to the determinant $|\L|$;
$\si_\ell$
takes $[a]$ to $[\ell a]$, and all parities $\eps_\ell([a])=+1$.
The Galois action for the affine algebras is quite interesting (see e.g.\
[1]), and can be expressed
geometrically using the action of the affine Weyl group on the weight lattice
of $X_r$. Both $\eps_\ell(\la)=\pm 1$ will occur. For finite groups,
$\si_\ell$ takes $(a,\chi)$ to $(a^\ell,\si_\ell\circ\chi)$, and again all
$\eps_\ell(a, \chi)=+1$.

The presence of the Galois action (4.4) is an effective criterion (necessary
and sufficient) on the matrix $S$ for the numbers in (2.1) to be
rational. It would be very desirable to find effective conditions on  $S$
such that the fusion coefficients are nonnegative, or integral. At present
the best
results along these lines are, respectively, the inequalities (4.2), and
the fact that the ratios ${S_{ab}\over S_{0b}}$ are algebraic integers
(since they are eigenvalues of integer matrices). When there are units,
then (4.3a) provides an additional strong constraint on nonnegativity.

Whenever a structure is studied, of fundamental importance are the
structure-preserving maps. It is through these maps that different
examples of the structure can be compared.
By a {\it fusion-homomorphism} $\pi$ between fusion rings
${\cal  F}(\Phi,N)$ and ${\cal F}(\Phi',N')$
we mean a
ring homomorphism for which $\pi(\Phi)\subseteq\Phi'$. 
{\it Fusion-isomorphisms} and {\it fusion-automorphisms} are defined in the obvious
ways. All fusion-isomorphisms between affine algebra fusion rings are known.
Most of them are in fact fusion-automorphisms, and
are constructed in simple ways from the symmetries of the Dynkin
diagrams.  Here are some basic general facts about fusion-homomorphisms:

\medskip\noindent{{\smcap Proposition 3.}} {\it Let $\pi:\Phi\rightarrow\Phi'$ 
be a fusion-homomorphism between any two fusion rings. Then

\item{(a)} $\pi 0=0'$ and $\pi(a^*)=\pi(a)^*$, and $\pi$ takes units
of $\Phi$ to units of $\Phi'$.

\item{(b)} There exists a map $\pi':\Phi'\rightarrow\Phi$ such that
$${S'_{\pi a,b'}\over S'_{0',b'}}={S_{a,\pi'b'}\over S_{0,\pi'b'}}\qquad
\forall a\in\Phi,\ b'\in\Phi'$$

\item{(c)} If $\pi a=\pi b$, then $b=Ja$ for some simple-current $J$.
In addition, this $J$ will obey $\pi(Jd)=\pi(d)$ for all $d\in\Phi$, and
(provided $J$ is nontrivial) there can be no $J$-fixed-points in $\Phi$.

\item{(d)} If $\pi$ is surjective, then $\pi':\Phi'\rightarrow\Phi$ is
an injective fusion-homomorphism, and
$$S'_{\pi a,b'}=\sqrt{{\rm ker}(\pi)}\,S_{a,\pi'b'}$$

}

Part (a) follows from {\bf F1} and {\bf F3}. Part (b) follows because
${S'_{\pi a,b'}\over S'_{0',b'}}$ is a 1-dimensional representation
of the $\Phi'$ fusion ring. To get (c), consider $(\pi a)(\pi b)^*=
\pi(ab^*)$. If $f$ is a fixed-point of $J$ in (c),
count the multiplicity of the identity  $0'$ in the fusion
product $(\pi f)\cdot (\pi f)^*$. To see (d), apply (c) to 
$$\sum_a\left|{S'_{\pi a,b'}\over S'_{0',b'}}\right|^2=\sum_a\left|{S_{a,\pi'b'}\over
S_{0,\pi' b'}}\right|^2$$

For example, fix any units $j,j'\in \Phi$ of equal order $n$. 
Then $a\mapsto J^{Q'(a)}a$ defines a fusion-endomorphism, where we write
$\varphi_{j'}(a)=\exp[2\pi\i\, Q'(a)/n]$.  It will be a fusion-automorphism
iff $Q'(j)+1$ is coprime to $n$. For another example, take any Galois
automorphism $\si$ for which $\si(S_{00}^2)=S_{00}^2$, or equivalently
$\si 0=J0$ for some simple-current $J$. Then $a\mapsto J\si a$ is a
fusion-automorphism. For this Galois example $\pi'=\pi$, while for the
simple-current one $\pi'(b)=J'{}^{Q(b)}b$. 

The map $\pi'$ of Prop.\ 3(b) won't in general be a fusion-homomorphism. E.g.\ 
consider the fusion-homomorphism
$\pi:\{[0],[1]\}\rightarrow\{0,1,\ldots,k\}$ between the fusion ring of
the lattice $\L=\sqrt{2}\Z$ and the fusion ring for $A_1^{(1)}$ level
$k$, given by $\pi([0])=0,\pi([1])=k$. Then $\pi'$ is given by
$\pi'(a)= [a]$.

A very desirable property for modular data to possess is:

\medskip\noindent{\smcap Congruence Subgroup Property.} [21] {\it Let $N$ be the order of the matrix
$T$, so $T^N=I$, and let $\rho$ be the representation of SL$_2(\Z)$
coming from the assignment (1.2b). Then $\rho$ factors through the congruence subgroup
$$\Gamma(N):=
\{A\in {\rm SL}_2(\Z)\,|\,A\equiv\left(\matrix{1&0\cr 0&1}\right)\ ({\rm mod}\
N)\}$$
 and so (1.2b) in fact defines a representation of the finite group
SL$_2(\Z/N\Z)$. Moreover, the characters (1.1) are modular functions
for $\Gamma(N)$. The entries $S_{ab}$ all lie in the  cyclotomic field
$\Q[\exp(2\pi\i/N)]$, and for any Galois automorphism $\si_\ell$,}
$$T_{\si_\ell a,\si_\ell a}=T_{aa}^{\ell^2}\qquad \forall a\in\Phi\eqno(4.5)$$

For example, the modular data from Examples 1--3 in Section III all obey
this property. In particular, affine algebra characters $\chi_\la$
are essentially
lattice theta functions. It would be valuable to find examples of modular
data which do {\it not} obey this property. For much more discussion,
see [21,7].
In those papers, considerable progress was made towards clarifying its role (and
existence) in modular data. For example:

\medskip\noindent{\smcap Proposition 4.} [21] {\it Consider any modular data.
Let $N$ be the order of $T$, and suppose that $N$ is either coprime to
$p=2$ or $p=3$. Then the corresponding SL$_2(\Z)$ representation
factors through $\Gamma(N)$, provided $(4.5)$ holds for $\ell=p$.}\medskip

In the remaining case, i.e.\ when 6 divides $N$, more conditions are
needed; these are also given in [21].
It is tempting to think that this is a good approach to verifying that
rational VOA characters are modular functions. It also leads, via [32], to a promising approach to
classifying modular data.

Assuming some additional structure from RCFT,
[7] recently established the congruence property ([21] had previously
proved the $\Gamma(N)$ part when $T$ has odd order). Though this is clearly
an impressive feat, what it means in the more general context of modular data isn't
clear: it is difficult to explicitly write down
the additional axioms needed to supplement our definition of modular
data, in order that the necessary calculations go through.

\bigskip
\centerline{{\bf V. Modular Invariants and NIM-reps}}\medskip

A {\it modular invariant} is a matrix $M$, rows and columns labeled by
$\Phi$, obeying:

\smallskip\item{{\bf MI1.}} $MS=SM$ and $MT=TM$;

\smallskip\item{{\bf MI2.}} $M_{ab}\in\Z_{\ge}$ for all $a,b\in\Phi$; and

\smallskip\item{{\bf MI3.}} $M_{00}=1$.\smallskip

As usual we write $\Z_\ge$ for the nonnegative integers. The simplest
example of a modular invariant is of course the identity matrix $M=I$. Another
example is conjugation $C$. All of the modular invariants for $A_1^{(1)}$ at level
$k$ are given below in (6.1).

Why are modular invariants interesting? Most importantly, they are central to
the task of classifying RCFTs. The genus-1 `vacuum-to-vacuum amplitude'(=partition
function) ${\cal Z}(\tau)$ of the theory looks like (1.3c). It assigns to
the torus $\C/(\Z+\Z\tau)$ the complex number ${\cal Z}(\tau)$.
But the moduli space of conformally equivalent tori is the orbit space
SL$_2(\Z)\backslash{\Bbb H}$,
where the action is given by $\left(\matrix{a&b\cr c&d}\right)\tau=
{a\tau+b\over c\tau+d}$. Thus the partition function
${\cal Z}(\tau)$ must be invariant under this natural action of the modular group
SL$_2(\Z)$, which gives us {\bf MI1}.
The coefficients $M_{ab}$ count the primary fields
$|\phi_a,\phi_b\rangle$ in the state space ${\cal H}$, i.e.\ the number of times the
module $A_a\otimes A_b$ of left chiral algebra$\times$right chiral
algebra, appears in ${\cal H}$. That gives 
us {\bf MI2}. And the uniqueness of the vacuum $|0,0\rangle$ means {\bf MI3}. That is to say,
the coefficient matrix $M$ of an RCFT partition function is a modular
invariant. It is believed that an RCFT is uniquely specified by the knowledge
of its partition function, its (left and right) chiral algebras($=$VOAs), and
the so-called structure constants. In any case, an important fingerprint
of the RCFT is its partition function ${\cal Z}$, i.e.\ its
modular invariant $M$.

Another motivation for studying modular invariants is the extensions
${\cal V}\subset {\cal V}'$ of rational VOAs (similar remarks hold for
braided subfactors). Let $M_i$ and $M'_j$ be the irreducible modules
of ${\cal V}$ and ${\cal V}'$, respectively. Then each $M'_j$ will be
a ${\cal V}$-module. A rational
VOA should have the complete reducibility property, so each $M'_j$ should
be expressible as a direct sum of $M_i$'s --- these are called the branching
rules. As mentioned in Example 5, we would expect that the characters (1.1) of a rational VOA should yield
(some form of) modular data via (1.2a). So the diagonal sum $\sum_j|{\rm ch}'_{M'_j}|^2$
should be invariant under the SL$_2(\Z)$-action; rewriting the ${\rm ch}'_{M'_j}$'s
in terms of the ${\rm ch}_{M_i}$'s via the branching rules yields a modular invariant
for ${\cal V}$.

For instance, the VOA $L(\L_0)'$ corresponding to the affine algebra
$G_2^{(1)}$ level 1 contains the VOA $L(28\L_0)$ corresponding to $A_1^{(1)}$
at level 28. We get the branching rules
$L(\L_0)'=L(0)\oplus L(10)\oplus L(18)\oplus L( 28)$ and
$L(\L_2)'=L(6)\oplus L( 12)\oplus L(16)\oplus L(22)$, where $L(\la_1):=L(\la)$. This corresponds to
the $A_1^{(1)}$ level 28 modular invariant given below in (6.1f).

So knowing the modular invariants for some VOA ${\cal V}$ gives considerable information
concerning its possible `nice' extensions ${\cal V}'$. For instance, we are learning that the only finite
`rational' extensions of a generic  affine VOA are those studied in [29] (`simple-current
extensions') and whose modular data is conjecturally given in [41].

Another reason for studying modular invariants is that the answers are
often surprising. Lists arising in math from complete classifications tend to
be about as stale as phonebooks, but to give some samples:

\smallskip\item{$\bullet$} the $A_1^{(1)}$ modular invariants fall into the A-D-E
metapattern;

\item{$\bullet$} the $A_2^{(1)}$ modular invariants have connections with
Jacobians of Fermat curves; and

\item{$\bullet$} the $(U(1)\oplus\cdots\oplus U(1))^{(1)}$
modular invariants correspond to rational points on Grassmannians.

\smallskip\noindent We will
discuss this point a little more next section. 
 These `coincidences', presumably, have something to do with the nontrivial
 connections between RCFT and several areas of math, but it also is due to
the beauty of the combinatorics of Lie characters evaluated at elements of finite
order (3.2c).

In any case, in this section we will  study
the modular invariants corresponding to a given choice of modular data.
For lattices, the classification is easy (use (5.2) below). For many finite groups,
the classification typically will be hopeless --- e.g.\ the alternating group
$A_5$, which has only 22 primaries, has a remarkably high number (8719) of
modular invariants [5].
For affine algebra modular data, the classification of modular invariants seems
to be just barely possible, and the answer is that (generically) the only
modular invariants are constructed in straightforward ways from symmetries
of the Dynkin diagrams.\smallskip

Commutation with $T$ is trivial to solve, since $T$ is diagonal: it
yields the selection rule
$$M_{ab}\ne 0\ \Rightarrow\ T_{aa}=T_{bb}\eqno(5.1)$$
This isn't as useful as it looks; commutation with $S$ (or
equivalently, the equation $SM\overline{S}=M$) is more subtle,
but far more valuable. 

An immediate observation is that there are only finitely many modular invariants
associated to given modular data. This follows for instance from
$$1=M_{00}=\sum_{a,b\in\Phi}S_{0a}\,M_{ab}\,S_{b 0}
\ge S_{00}^2\sum_{a,b\in\Phi}M_{ab}$$

We will find that each basic symmetry of
the $S$ matrix yields a symmetry of the modular invariants, a selection
rule telling us that certain entries of $M$ must vanish, and a way to
construct new modular invariants.

First consider simple-currents $J,J'$. Equation (4.3a) and positivity tell us
$$M_{J0,J'0}=\left|\sum_{c,d\in\Phi}\varphi_J(c)\,S_{0c}\,M_{cd}\,\overline{S_{d0}}
\,\overline{\varphi_{J'}(d)}\right|\le \sum_{c,d}S_{0c}M_{cd}S_{d0}=M_{00}=1\eqno(5.2a)$$
Thus $M_{J0,J'0}\ne 0$ implies $M_{J0,J'0}=1$, as well as the selection rule
$$M_{cd}\ne 0\ \Rightarrow\ \varphi_J(c)=\varphi_{J'}(d)\eqno(5.2b)$$
A similar calculation yields the symmetry
$$M_{J0,J'0}\ne 0\ \Rightarrow\ M_{Ja,J'b}=M_{ab}\qquad \forall a,b\in\Phi\eqno(5.2c)$$

The most useful application of simple-currents to modular invariants
is to their construction. In particular, let $J$ be a simple-current of
order $n$. Then we learned in (4.3) that $\varphi_j(a)$ is an $n$th
root of 1, and that
$(T_{jj}\overline{T_{00}})^{2n}=1$ and in fact
$(T_{jj}\overline{T_{00}})^{n}=1$ when
$n$ is odd. That is to say,  we can find integers  $r_j$ and $Q_j(a)$ obeying
$$\vp_j(a)=\exp[2\pi\i\,{Q_j(a)\over n}]\ ,\qquad T_{jj}\,\overline{T_{00}}=\exp[\pi\i \,r_j{n-1\over n}]$$
For $n$ odd, choose $r_j$ to be even (by adding $n$ to it if necessary).
Now define the matrix ${\cal M}[J]$ by [81]
$${\cal M}[J]_{ab}=\sum_{\ell=1}^n \delta_{J^\ell a,b}\,\delta({Q_j(a)\over n}+
{\ell\over 2n}r_j)\eqno(5.3)$$
where $\delta(x)=1$ when $x\in\Z$ and is 0 otherwise. This matrix
${\cal M}[J]$
will be a modular invariant iff $(T_{jj}\overline{T_{00}})^n=1$ (i.e.\
iff $r_j$ is even), and a permutation matrix iff $T_{jj}\overline{T_{00}}$
is a {\it primitive} $n$th root of 1.
When $n$ is even, (4.3c) says $(T_{jj}\overline{T_{00}})^n=1$  iff
$\varphi_j(j)^{n/2}=1$.

For instance, taking
$J=id$ we get ${\cal M}[id]=I$. The affine algebra $A_1^{(1)}$ at level $k$
has a simple-current  with $r_j=k$ given by $Ja=k-a$; for even $k$ the
matrix ${\cal M}[J]$ is 
the modular invariant called ${\cal D}_{{k\over 2}+2}$ below in (6.1b),(6.1c).

Now look at the consequences of Galois. Applying the Galois automorphism
$\si$ to $M=SM\overline{S}$ yields from (4.4) and $M_{ab}\in\Q$ the fundamental
equation
$$M_{ab}=\sum_{c,d\in\Phi}\eps_\si(a)\,S_{\si a,c}\,M_{cd}\,\overline{S_{d,\si b}}
\,\eps_\si(b)=\eps_\si(a)\,\eps_\si(b)\,M_{\si a,\si b}\eqno(5.4a)$$
Because $M_{ab}\ge 0$, we obtain the selection rule
$$M_{ab}\ne 0\ \Rightarrow\ \eps_\si(a)=\eps_\si(b)\qquad \forall\si\eqno(5.4b)$$
and the symmetry
$$M_{\si a,\si b}=M_{ab}\qquad \forall\si\eqno(5.4c)$$
Of all the equations (5.2) and (5.4), (5.4b) is the most valuable. A way to construct
modular invariants from Galois was first given in [40] but isn't useful
for constructing affine algebra modular invariants and so won't be repeated here.

There are other very useful facts, which space prevents us from describing.
For instance, we have the inequality
$$\sum_{b\in\Phi}S_{ab}M_{b0}\ge 0\eqno(5.5)$$
Perron-Frobenius tells us many things, e.g.\ that any modular invariant $M$
obeying $M_{0a}=\delta_{0a}$ must be a permutation matrix. For
affine algebra modular invariants, the Lie theory of the underlying finite-dimensional
Lie algebra plays a crucial role, thanks largely to (3.2c).

Closely related to modular invariants is the notion of {\it NIM-rep} (short
for `nonnegative integer representation' [12]) or equivalently
{\it fusion graph}. These originally arose in two {\it a priori} unrelated contexts:
the analysis, starting with Cardy's fundamental paper [16], of boundary RCFT;
and Di Francesco--Zuber's largely empirical attempt
[24] to understand and generalise the A-D-E metapattern appearing in $A^{(1)}$
modular invariants, by attaching graphs to each conformal field theory.

A {\it NIM-rep} $\N$ is a nonnegative integer representation of the fusion ring,
that is, an assignment $a\mapsto \N_a$ to each $a\in\Phi$ of a matrix $\N_a$
with nonnegative integer entries, obeying $\N_a\N_b=\sum_c N_{ab}^c\N_c$.
In addition we require that $\N_0=I$ and that transpose and conjugation
be related by $\N_a^t=\N_{Ca}$, for all $a\in\Phi$.

Two obvious examples of NIM-reps are the fusion matrices, $a\mapsto N_a$,
and their transposes $a\mapsto N_a^t$. The rows and columns of most NIM-reps
however won't be labelled by $\Phi$, in fact we will see that the dimension
of the NIM-rep should equal the trace Tr$(M)$ of some modular invariant.

Just as it is convenient to replace a Cartan matrix by its Dynkin diagram, so
too is it convenient to realise $\N_a$ by a (directed multi-)graph: we put a
node for each row/column, and draw $(\N_a)_{\alpha\beta}$ edges directed
from $\alpha$ to $\beta$. We replace each pair of arrows $\alpha\rightarrow
\beta,\beta\rightarrow \alpha$, with a single undirected edge connecting
$\alpha$ and $\beta$. These graphs are called {\it fusion
graphs}, and are often quite striking.

NIM-reps correspond in RCFT  to the  1-loop vacuum-to-vacuum 
amplitude ${\cal Z}_{\alpha\beta}(t)$ of an open string, or
equivalently 
of a cylinder whose edge circles are labelled by `conformally invariant boundary
states' $|\alpha\rangle,|\beta\rangle$ [16,80,42,11]. In string theory these
are called the `Chan-Paton degrees-of-freedom' and are placed at the endpoints
of open strings. The real variable $-\infty<t<\infty$ here is the modular parameter
for the cylinder, and plays the same role here that $\tau\in{\Bbb H}$ plays in
 ${\cal Z}(\tau)$. In particular we get (1.4), where the matrices
 $(\N_a)_{\alpha\beta}=\N_{a\alpha}^\beta$ define a NIM-rep.
These (finitely many) boundary states $\alpha$ are the indices for
the rows and columns of each matrix $\N_a$.

By the usual arguments (see Section IV) we can simultaneously diagonalise all
$\N_a$, and the eigenvalues of $\N_a$ will be $S_{ab}/S_{0b}$ for $b$ in
some multi-set $\E=\E(\N)$ (i.e.\ the elements of $\E$ come with
multiplicities). This multi-set $\E$ depends only on $\N$ (i.e.\ is
independent of $a\in\Phi$), and is called the {\it exponents} of the NIM-rep.

Two NIM-reps $\N,\N'$ are regarded as equivalent if there is a simultaneous
permutation $\pi$ of the rows and columns such that $\pi\N_a\pi^{-1}=\N_a'$
for all $a\in\Phi$. For example, the two NIM-reps given earlier are equivalent:
$N_a^t=CN_aC^{-1}$. We write $\N=\N'\oplus\N''$, and call $\N$ {\it reducible},
 if the matrices $\N_a$ can be simultaneously written as direct sums
$\N_a=\N_a'\oplus\N''_a$. Necessarily, the summands
 $\N'$ and $\N''$ themselves will  be  NIM-reps. 
Irreducibility is equivalent to demanding that the identity 0 occurs in $\E(\N)$
with multiplicity 1. We are interested in {irreducible} equivalence 
classes of NIM-reps --- there will be only finitely many [53].

Two useful facts are: the Perron-Frobenius
eigenvalue of $\N_a$ is the q-dimension ${S_{a0}\over S_{00}}$ (we'll see
this used next section); and for all $a\in\Phi$,
$$\sum_{b\in{\cal E}}{S_{ab}\over S_{0b}}={\rm
Tr}\,(\N_a)\in\Z_{\ge}\eqno(5.6) $$
The consequences of the simple-current and Galois symmetries are also
important and are worked out in [53].

By the {\it exponents} of a modular invariant $M$ we mean the multi-set $\E_M$
where $a\in\Phi$ appears with multiplicity $M_{aa}$. RCFT [16,11] is thought to
require that each modular invariant $M$ have a companion NIM-rep
$\N$ with the property that 
$$\E_M=\E(\N) \eqno(5.7)$$ 
So the size of the matrices
$\N_a$, i.e.\ the dimension of the NIM-rep, should equal the trace Tr$(M)$
of the modular invariant. For instance,
the fusion matrix NIM-rep $a\mapsto N_a$ corresponds to the modular invariant
$M=I$. However, there doesn't seem to be a general expression for the
NIM-rep (if it exists) of the next simplest modular invariant, the
conjugation $M=C$. 

Incidentally, the inequality (5.6) is
automatically obeyed by the exponents ${\cal E}=\E_M$ of any modular
invariant $M$:
$$\sum_{b\in\E_M}{S_{ab}\over S_{0b}}={\rm Tr}(MD_a)={\rm Tr}(\overline{S}SMD_a)=
{\rm Tr}(MSD_a\overline{S})={\rm Tr}(MN_a)\in\Z_\ge$$

Note that the NIM-rep definition depends on $S$, while a modular
invariant also sees $T$. One consequence of this is the
following. Suppose there is a primary $a\in\Phi$ such that 
$$T_{bb}=T_{cc}\ \Rightarrow\ S_{ab}\,\overline{S_{ac}}\ge 0\qquad
\forall b,c\in\Phi\eqno(5.8)$$
Then $M_{aa}=\sum_{b,c}S_{ab}M_{bc}\overline{S_{ac}}>0$ and so
$a\in{\cal E}_M$. It is thus natural to require of a NIM-rep $\N$ that
any such primary $a\in\Phi$ must appear in ${\cal E}(\N)$ with
multiplicity $\ge 1$, because otherwise no modular invariant $M$ could
be found obeying (5.7).

It should be mentioned that, from the RCFT point of view, the
constraint (5.7) is not as `carved  in stone' as {\bf MI1}--{\bf MI3}.
Our treatment here of NIM-reps reflects the current understanding, but it is still
based on unproven physical assumptions (`completeness of boundary
conditions') and perhaps in the future will
require some modification. Also, we're ignoring here the `pairing'=`gluing
automorphism' $\omega$ of e.g.\ [42], although here this isn't a serious omission.
 But an independent justification for studying
NIM-reps, and a strong hint that this RCFT picture is not too naive, comes
from subfactors.

NIM-reps and modular invariants appear very naturally in the  subfactor picture (Example 6)
[33,76,12], again paired by the relation (5.7). 
In this remarkable picture it is possible to interpret not only
the diagonal entries of the modular invariant, but in fact
all entries [75,12] (this was already anticipated in [24]). Extend the setting of Example 6 by considering a
braided system of endomorphisms for  a type III subfactor $N\subset M$.
Here, the primaries $\Phi={}_N\chi_N$ consist of irreducible
endomorphisms of $N$, while the rows and columns of our NIM-rep will
be indexed by irreducible homomorphisms $a\in{}_M\chi_N$, $a:N\rightarrow
M$. The fusion-like ring of ${}_N\chi_N$ will be commutative, i.e.\ be a true
fusion ring; that of ${}_M\chi_M$ however will generally be
noncommutative. There is a simple expression [12] for the corresponding
modular invariant using `$\alpha$-induction' (a process of inducing an
endomorphism from $N$ to $M$ using the braiding $\epsilon^{\pm}$):
we get $M_{\la\mu}=\langle \alpha_\la^+,\alpha_\mu^-\rangle$ where the
dimension $\langle,\rangle$ is defined in Example 6. Then the
(complexified) fusion algebra of ${}_M\chi_M$ will be isomorphic (as a
complex algebra) to $\oplus_{\la,\mu}{\rm GL}_{M_{\la\mu}}({\Bbb
C})$.  The NIM-rep is essentially $\alpha$:
 $({\cal N}_\la)_{a,b}=\langle b,\alpha^{\pm}_\la a\rangle$ (either
choice of $\alpha^\pm$ gives the same matrix) [12].  This NIM-rep
arises as a natural action of ${}_M\chi_M$ on ${}_M\chi_N$.
As these partition functions of tori and cylinders appear so nicely
here, it is tempting to ask about other surfaces, especially the
M\"obius band and Klein bottle, which also play a basic role in
boundary RCFT [80].

We won't speak much more here about NIM-reps --- see e.g.\ [24,11,53] and references
therein for more of the theory and classifications (and graphs!). Typically, what has happened in the classifications
thus far is that there are slightly
more NIM-reps than modular invariants, but their classifications match
surprisingly well. For instance, the irreducible NIM-reps of $A_1^{(1)}$
have $\N_1$ equal to the incidence matrix of the A-D-E graphs and
tadpoles [24]
--- compare with the list of modular invariants for $A_1^{(1)}$ in (6.1)! However
there are places where many modular invariants lack a corresponding NIM-rep
(this happens for instance for the orthogonal algebras at level 2 [53]).
The simplest examples of modular invariants lacking NIM-reps
occur for $B_4^{(1)}$ level 2, and the symmetric group $S_3$.

A tempting guess is that almost all of the enormous numbers of modular
invariants associated to finite group modular data will likewise fail to have a
corresponding NIM-rep. Recall that the Galois parities $\epsilon_\ell$
for the finite group modular data are all $+1$, and hence the
constraint (5.4b) becomes trivial. As a general rule, the number of
modular invariants  is inversely
related to the severity (5.4b) possesses for that choice of modular data.

The moral of the story seems to be the following. The definition of modular
invariants didn't come to us from God; it came to us from men like Witten, Cardy, ...
The surprising thing is that so often their classification yields interesting
answers. A modular invariant
may not correspond to a CFT (we have infinitely many examples where it fails
to), and the modular invariant may correspond to different CFTs
(though all known examples of this are artificial, due to our characters
depending on too few variables to distinguish the representations of the
maximally extended VOAs). But ---
at least for most affine algebras and levels --- it seems they're {\it usually} in
one-to-one correspondence.

In any case, classifying modular invariants, and comparing their lists
to those of NIM-reps, is a natural task and has led to interesting findings
(see e.g.\ the review [93]). 

\bigskip\centerline{{\bf VI. Affine Algebra Modular Invariant Classifications}}\medskip

The most famous modular invariant classification was the first. In (3.5)
we gave explicitly the modular data for the affine algebra $A_1^{(1)}$
at level $k$.
Its complete list of modular invariants is [15] (using the simple-current $Ja=k-a$)
$$\eqalignno{{\cal A}_{k+1}=&\,\sum_{a=0}^k|\chi_a|^2\ ,\qquad \forall k\ge 1&(6.1a)\cr
{\cal D}_{{k\over 2}+2}=&\,\sum_{a=0}^k\chi_a\,\overline{\chi_{J^aa}}\ ,\qquad
{\rm whenever}\ {k\over 2}\ {\rm is\ odd} &(6.1b)\cr
{\cal D}_{{k\over 2}+2}=&\,|\chi_0+\chi_{J0}|^2+|\chi_2+\chi_{J2}|^2+\cdots
+2|\chi_{{k\over 2}}|^2\ ,\qquad{\rm whenever}\ {k\over 2}\ {\rm is\ even}
&(6.1c)\cr
{\cal E}_6=&\,|\chi_0+\chi_6|^2+|\chi_3+\chi_7|^2+|\chi_4+\chi_{10}|^2\ ,\qquad
{\rm for}\ k=10&(6.1d)\cr
{\cal E}_7=&\,|\chi_0+\chi_{16}|^2+|\chi_4+\chi_{12}|^2+|\chi_6+\chi_{10}|^2
&\cr&\,+\chi_8\,(\overline{\chi_2+\chi_{14}})+(\chi_2+\chi_{14})\,
\overline{\chi_8}+|\chi_8|^2\ ,\qquad {\rm for}\ k=16&(6.1e)\cr
{\cal E}_8=&\,|\chi_0+\chi_{10}+\chi_{18}+\chi_{28}|^2+|\chi_{6}+\chi_{12}
+\chi_{16}+\chi_{22}|^2\ ,\ {\rm for}\ k=28&(6.1f)}$$
Each of these is identified with a (finite) Dynkin diagram, in such a way that
the Coxeter number $h$ of the diagram equals $k+2$, and the {\it exponents}
of the corresponding Lie algebra are given by $1+{\cal E}_M$ (recall the
definition of exponents ${\cal E}_M$ of a modular invariant, given at the
end of last section). The exponents of the Lie algebra are the numbers $m_i$,
where $4\sin^2(\pi\,{m_i\over h})$ are the eigenvalues of the Cartan
matrix. For instance, the Dynkin diagram $D_8$ has Coxeter number 14 and exponents
$1,3,5,7,7$, while ${\cal D}_8$ occurs at level 12 and has exponents
${\cal E}=\{0,2,4,6,6\}$.

The A-D-E pattern appears in many places in math and mathematical physics [60]:
besides the simply-laced Lie algebras and $A_1^{(1)}$ modular invariants,
these diagrams also classify simple singularities, finite subgroups
of SU$_2({\Bbb C})$, subfactors with Jones index $<4$, representations
of quivers, etc. There seem to be two more-or-less inequivalent A-D-E patterns,
one corresponding to the finite A-D-E diagrams, and the other corresponding to
the affine (=extended) A-D-E diagrams. For instance, the modular invariants
identify with the finite ones, while the finite subgroups of SU$_2({\Bbb C})$
match with the affine ones. This suggests that a direct relation between e.g.\
the modular invariants and those finite subgroups could be a little forced.
Patterns such as A-D-E are usually explained by identifying an underlying
combinatorial fact which is responsible for its various incarnations.
The A-D-E combinatorial fact is probably the classification of
symmetric matrices over $\Z_{\ge}$, with no diagonal entries, and
 with maximal eigenvalue $<2$ (for the finite diagrams) and $=2$
(for the affine ones). Perhaps the only A-D-E classification which
still resists this `explanation' is that of $A_1^{(1)}$ modular invariants.
This is in spite of considerable effort (and some progress) by many people.
The present state of affairs, and also a much simpler proof on the lines
sketched in the previous section, is provided by [51].

Many other classes of affine algebras and levels have been classified.
The main ones are: $A_2^{(1)}$, $(A_1+A_1)^{(1)}$, and $(U(1)+\cdots+U(1))^{(1)}$,
 for all levels $k$; and 
$A_r^{(1)},B_r^{(1)},D_r^{(1)}$ for all ranks $r$, but with levels restricted to
$k\le 3$. See e.g.\ [50] for references to these results.

Has A-D-E been spotted in these other lists? No. However, a remarkable
connection [79] has been observed between the $A_2^{(1)}$ level $k$
modular invariants,  and the Jacobian of the Fermat curve $x^{k+3}+y^{k+3}
+z^{k+3}=0$. In particular, the $A_2^{(1)}$ Galois selection rule (5.4b)
and the analysis of the  simple factors in the Jacobian are essentially
the same. This link between Fermat and $A_2^{(1)}$ is still unexplained, and
how it extends to the other algebras,
e.g.\ perhaps $A_r^{(1)}$ level $k$ relates to $x_1^{k+r+1}+x_2^{k+r+1}+\cdots+x_{r+1}^{k+r+1}=0$?,
is still unclear. However, Batyrev [8] has suggested some possibilities
involving toric geometry.

The third `sample' listed last section (relating $(U(1)\oplus\cdots\oplus
U(1))^{(1)}$ modular invariants to the Grassmannians) suggests a different
link with geometry. The Grassmannian is essentially the moduli space of
Narain compactifications of a (classical) string theory, so perhaps
other families of modular invariants can be regarded as special
points on other finite-dimensional moduli spaces.

Though there are no other appearances of A-D-E, there is  a rather natural way to assign (multi-di)graphs to modular
invariants, generalising the A-D-E pattern for $A_1^{(1)}$. Note first that
we can classify the $A_1^{(1)}$ NIM-reps [24]: ${\cal N}_1$ must be symmetric
and have Perron-Frobenius eigenvalue ${S_{10}\over S_{00}}=2\cos({\pi\over k+2})<2$; thus the graph
associated to $\N_1$ must be an A-D-E Dynkin diagram, or a tadpole. The
tadpoles can be discarded, since they don't correspond via (5.7) to a modular invariant.
Given $\N_1$, all other $\N_a$ can be recursively obtained using the
special case $\N_1\N_i
=\N_{i+1}+\N_{i-1}$ of (3.5c). The result is a NIM-rep.

In this way, we find that the Dynkin diagram which (6.1) assigned to a
given
$A_1^{(1)}$ modular invariant $M$ is precisely the graph whose adjacency matrix
equals the generator $\N_1$ of the unique NIM-rep compatible with $M$
in the sense of (5.7).
Likewise, we should assign to the modular invariants of e.g.\ $A_2^{(1)}$
the multi-digraph $\N_{\L_1}$ generating the corresponding NIM-rep. The
NIM-reps for $A_2^{(1)}$ are not yet classified, but at least one has been
found for each $M$ [24,76,11,12].

There is a simple reason why the tadpole can't correspond to an
$A_1^{(1)}$ modular invariant. Note that the unit $a=k$ satisfies
(5.8), and thus will lie in any ${\cal E}_M$. However, $k$ is not an
exponent of the tadpole, and thus there can be no solution $M$ in
(5.7) for the choice $\N=\,$tadpole. More generally, this suggests
refining the definition of NIM-rep: many extraneous (unphysical?)
NIM-reps can be avoided, by requiring $a\in{\cal E}(\N)$ for any
$a\in\Phi$ satisfying (5.8).

By the way, {\it sub}modular invariants can usually be found for
NIM-reps which lack a true modular invariant. For example, the
seemingly extraneous $n$-vertex tadpole mentioned in the previous
paragraph corresponds to the algebra $A_1^{(1)}$ at level $2n-1$, and
the submodular invariant $M_{ab}=\delta_{b,J^aa}$.  Perhaps a reasonable
interpretation can be found by both the subfactor and boundary CFT camps
for  NIM-reps corresponding to matrices $M$ commuting with certain small-index subgroups of SL$_2(\Z)$.
 Recall that we anticipated this thought at the end of Example 1.

Most of the modular invariant classification effort has been directed not at specific algebras and levels, but at the
general argument. The major result obtained thus far is:

\medskip\noindent{{\smcap Theorem 5.}} [49] {\it Choose any affine
algebra $X_r^{(1)}$ and level $k$. Let $M$ be any modular invariant, 
 obeying the constraint that the only primaries $a\in\Phi$ for which
$M_{0a}\ne 0$ or $M_{a0}\ne 0$, are units. Then $M$ lies on an explicit list.}
\medskip

Note that, of the $A_1^{(1)}$ modular invariants, all but
 ${\cal E}_6$ and ${\cal E}_8$ obey the constraint of Thm.\ 5. That pattern
 seems to continue for the other algebras and levels: the list of modular
 invariants covered by Thm.\ 5 exhausts almost every modular invariant
 yet discovered.

There are very few {\it exceptional} modular invariants in the list of Thm.\ 5.
Almost all of the modular invariants there are simple-current ones (5.3), and the product
of these by the conjugation $C$ (strictly speaking, any symmetry of
the unextended Dynkin diagram can be used here in place of $C$). 

Thm.\ 5 is important because, {\it for generic choice of algebra and level,} the
various constraints we have on the 0-row and 0-column of a modular invariant
(most importantly, Galois (5.4b), $T$ (5.1), and the inequality (5.5))
force the condition of Thm.\ 5 to be satisfied.

Indeed, if we impose the full structure of Ocneanu cells [76] (this should be
equivalent to saying that an RCFT exists with partition function given by $M$),
we obtain Ocneanu's inequality:
$$\sum_{\mu\in {\rm clearing}}N_{\la,C\la}^\mu\, S_{\mu 0}\le S_{\la 0}\eqno(6.2)$$
where $\la$ is any weight$\ne 0$ obeying $M_{\la 0}\ne 0$ with $\la_0$
as large as possible, and where
`clearing' is a subset of $P_+^k$ close to 0: $\mu$ is in the clearing
if $2(k-\mu_0)\le k-\la_0$. The left-side of (6.2) grows approximately
quadratically with $S_{\la 0}/S_{00}$, while the right-side is only linear,
so it tends to force $S_{\la 0}$ to be small; equation (5.1) on the
other hand tends to force $S_{\la 0}$ to be large.
This should imply that, for fixed algebra $X_r^{(1)}$, there is a $K$ (depending
on the algebra) such that $\forall k>K$, the constraint of Thm.\ 5 will
be obeyed! Thus:

\medskip\noindent{{\smcap Corollary 6.}}  {\it All possible modular invariants
appearing in RCFT (or the subfactor interpretation), corresponding to
any fixed choice of affine algebra $X_r^{(1)}$, and all sufficiently high
levels, are known.}\medskip

In other words, what Cor.\ 6 tells us is that, apart from some low level
exceptional modular invariants, all affine algebra modular invariants appearing in
RCFT can be constructed in straightforward and known ways from the
symmetries of the corresponding affine Dynkin diagram!

Thm.\ 5 has another consequence. It makes it relatively easy to find all
modular invariants (using only conditions {\bf MI1}-{\bf MI3})
at `small' levels, when the rank of the algebra isn't
too large [54]. For example, all modular invariants for $E_8^{(1)}$ at all
levels $k\le 380$ can be determined. This isn't completely trivial:
$E_8^{(1)}$ at $k=380$ has over $10^{12}$ highest weights$=$primaries,
so the $S$ and $M$ matrices have a number of entries approximately
equal to Avogadro's number! And each of these entries of $S$, given by
(3.2b), involves a sum of $10^9$ complex numbers. The fact that we can
reach such high levels isn't a sign of programming prowess, but rather
to how close we are to a complete classification of these (unrestricted)
affine algebra modular invariants. In [54] the modular invariants are given
for all exceptional algebras, and the classical algebras of rank $\le 6$.

The big surprise here is how rare the affine algebra modular invariants are
(for comparison, recall that there are over 8000 modular invariants for the
finite group $A_5$).
In the Table we've summarised the modular invariant classifications for
various algebras of small rank. It describes the complete list of
modular invariants for these algebras, when the level is sufficiently
small (these limits are given in the Table). 
A very safe conjecture though is that the Table gives the complete
classification for those algebras, for {\it all} levels $k$ (at the time of writing, $E_7^{(1)}$
level 42 and $E_8^{(1)}$ level 90 still have not been eliminated). 
Our hope is that this Table (or more realistically, the paper [54] where more
results are given and in more detail) will inspire
someone to spot a new coincidence involving modular invariants and some
other area of mathematics. For example, note in $A_1^{(1)}$ that the exceptionals
appear at $k+2=12,18,30$, which are the Coxeter numbers of $E_6$, $E_7$, $E_8$.
Claude Itzykson noticed that the $A_2^{(1)}$ exceptionals occur at
$k+3=8,12,24$ ---all divisors of 24--- and (inspired by the Fermat
connection [79])
found signs of these exceptionals in the Jacobian of $x^{24}+y^{24}+
z^{24}=0$. Can anyone spot any such pattern for the other algebras?

\bigskip\noindent{{\smcap Acknowledgements.}} I have benefitted from
several conversations with, and/or comments from, E.\ Bannai, D.\ Evans,
M.\ Gaberdiel, J.\ Lepowsky, V.\ Linek, P.\ Ruelle, M.\ Walton, and
J.-B.\ Zuber. This paper was written at St.\ John's
College in Cambridge, an ideal environment for work, and I thank them
and especially P.\ Goddard for generous hospitality.
 
\bigskip\centerline{{\bf References}}\medskip

\item{1.} D.\ Altschuler, P.\ Ruelle, and E.\ Thiran, ``On parity functions
in conformal field theories'', {\it J.\ Phys.\ A: Math.\ Gen.} {\bf 32} (1999),
3555--3570.

\item{2.} M.\ Atiyah, ``Topological quantum field theories'', {\it Publ.\
Math.\ IHES} {\bf 68} (1989), 175--186.

\item{3.} {H.\ Awata} and {Y.\ Yamada}, ``Fusion
rules for the fractional level $\widehat{sl(2)}$ algebra'', {\it Mod.\
Phys.\ Lett.} {\bf A7} (1992), 1185--1195;

\item{} {B.\ Feigin} and {F.\ Malikov}, ``Modular functor and representation
theory of $\widehat{sl}(2)$ at a rational level'', In: {\it Operads},
Contemp.\ Math.\ 202, Amer.\ Math.\ Soc., Providence, 1997, pp.357--405.

\item{4.} {E.\ Bannai}, ``Association schemes and fusion algebras (an
introduction)'', {\it J.\ Alg.\ Combin.} {\bf 2} (1993), 327--344.

\item{5.} E.\ Bannai, E.\ Bannai, O.\ Shimabukuro, and M.\ Tagami,
``Modular invariants of the modular data of finite groups'', preprint.

\item{6.} {E.\ Bannai} and {T.\ Ito}, {\it Algebraic Combinatorics I:
Association Schemes}, Benjamin-Cummings, Menlo Park California, 1984.

\item{7.} P.\ Bantay, ``The kernel of the modular representation and the
Galois action in RCFT'', preprint (math.QA/0102149).

\item{8.} {V.\ Batyrev}, ``Vertex algebras and mirror symmetry'',
({\it talk, University of Warwick, October 1999}).

\item{9.} {L.\ B\'egin, P.\ Mathieu,} and
{M.A.\ Walton}, ``$\widehat{su(3)}_k$ fusion coefficients'',
{\it Mod.\ Phys.\ Lett.} {\bf A7} (1992), 3255--3265.

\item{10.} {L.\ B\'egin, A.N.\ Kirillov, P.\ Mathieu,} and
{M.A.\ Walton}, ``Berenstein-Zelevinsky triangles, elementary
couplings, and fusion rules'', {\it Lett.\ Math.\ Phys.} {\bf 28}
(1993), 257--268.

\item{11.} {R.E.\ Behrend, P.A.\ Pearce, V.B.\ Petkova,} and {
J.-B.\ Zuber}, ``Boundary conditions in rational
conformal field theories'', {\it Nucl.\ Phys.} {\bf B579} (2000), 707--773.

\item{12.} {J.\ B\"ockenhauer} and {D.E.\ Evans}, ``Subfactors
and modular invariants'', (math.OA/0008056);

\item{} {J.\ B\"ockenhauer} and {D.E.\ Evans}, ``Modular invariants
from  subfactors'', (math.OA/0006114).

\item{13.} R.\ E.\ Borcherds, ``Vertex algebras, Kac-Moody algebras, and the
Monster'', {\it Proc.\ Natl.\ Acad.\ Sci.\ (USA)} {\bf 83} (1986), 3068--3071.

\item{14.} R.\ E.\ Borcherds, ``What is Moonshine?'', In: {\it Proc.\
Intern.\ Congr.\ Math.}, Berlin, 1998, pp.607--615.

\item{15.} {A.\ Cappelli, C.\ Itzykson,} and {J.-B.\ Zuber},
``The A-D-E classification of minimal and $A^{(1)}_1$ conformal
invariant theories'', {\it Commun.\ Math.\ Phys.} {\bf 113} (1987), 1--26.

\item{16.} {J.\ Cardy}, 
``Boundary conditions, fusion rules and the Verlinde formula'',
{\it Nucl.\ Phys.} {\bf B324} (1989), 581--596.

\item{17.} {M.\ Caselle} and {G.\ Ponzano},
  ``Analyticity, modular invariance and the classification of three
  operator fusion algebras'', {\it Phys.\ Lett.}\ {\bf B242} (1990), 52--58.

\item{18.} H.\ Cohn, {\it A Classical Invitation to Algebraic Numbers and
Class Fields}, Springer, New York, 1988.

\item{19.} {J.H.\ Conway} and {N.J.A.\ Sloane}, {\it Sphere
packings, lattices and groups}, 3rd edn, Springer, Berlin, 1999.

\item{20.} {A.\ Coste} and {T.\ Gannon}, ``Remarks on
Galois in rational conformal field theories'', {\it Phys.\ Lett.} {\bf B323} (1994),
316--321.

\item{21.} {A.\ Coste} and {T.\ Gannon}, ``Congruence subgroups
and conformal field theory'', preprint (math.QA/0002044).

\item{22.} {A.\ Coste, T.\ Gannon} and {P.\ Ruelle}, ``Finite group
modular data'', {\it Nucl.\ Phys.} {\bf B581} (2000), 679--717.

\item{23.} {L.\ Crane} and {D.N.\ Yetter},
``Deformations of (bi)tensor categories'', {\it Cah.\ Top.\ G\'eom.\
Diff.\ Cat\'eg.} {\bf 39} (1998), 163--180.

\item{24.} {P.\ Di Francesco} and {J.-B.\ Zuber}, ``SU($N$) lattice
integrable models associated with graphs'', {\it Nucl.\ Phys.} {\bf
B338} (1990), 602--646; 

\item{} {P.\ Di Francesco} and {J.-B.\ Zuber}, ``SU$(N)$
lattice integrable models and modular invariance'', In: {\it Recent
Developments in Conformal Field Theory} (1990), World Scientific, pp.179--215;

\item{25.} {Ph.\ Di Francesco, P.\ Mathieu,} and {D.\
S\'en\'echal}, {\it Conformal Field Theory}, Springer, New York, 1997.

\item{26.} {R.\ Dijkgraaf, C.\ Vafa, E.\ Verlinde} and
{H.\ Verlinde}, ``The operator algebra of orbifold models'',
{\it Commun.\ Math.\ Phys.} {\bf 123} (1989), 485--526.

\item{27.} {R.\ Dijkgraaf} and {E.\ Witten}, ``Topological gauge
theories
and group cohomology'', {\it Commun.\ Math.\ Phys.} {\bf 129} (1990), 393--429.

\item{28.} {C.\ Dong}, ``Vertex algebras associated with even lattices'',
{\it J.\ Alg.} {\bf 160} (1993), 245--265.

\item{29.} {C.\ Dong, H.\ Li,} and {G.\ Mason}, ``Simple currents
and extensions of vertex operator algebras'', {\it Commun.\ Math.\
Phys.} {\bf 180} (1996), 671--707.

\item{30.} {C.\ Dong, H.\ Li,} and {G.\ Mason}, ``Vertex
operator algebras associated to admissible representations of
$\widehat{sl_2}$'', {\it Commun.\ Math.\ Phys.} {\bf 184} (1997), 65--93.

\item{31.} {C.\ Dong, H.\ Li,} and {G.\ Mason}, ``Modular-invariance
of trace functions in orbifold theory and generalised moonshine'', {\it
Commun.\ Math.\ Phys.} {\bf 214} (2000), 1--56.

\item{32.} {W.\ Eholzer}, ``On the classification of modular
  fusion algebras'', {\it Commun.\ Math.\ Phys.}\ {\bf 172} (1995), 623--660.

\item{33.} {D.E.\ Evans} and {Y.\ Kawahigashi}, {\it Quantum
Symmetries on Operator Algebras}, Oxford University Press, 1998.

\item{34.} {G.\ Faltings}, ``A proof for the Verlinde formula'',
{\it J.\ Alg.\ Geom.} {\bf 3} (1994), 347--374.

\item{35.} A.J.\ Feingold and M.D.\ Weiner, ``Type A fusion rules from elementary
group theory'', preprint (math.QA/0012194).

\item{36.} {M.\ Finkelberg}, ``An equivalence of fusion categories'',
{\it Geom.\ Funct.\ Anal.} {\bf 6} (1996), 249--267.

\item{37.} {I.B.\ Frenkel, J.\ Lepowsky}, and {A.\ Meurman},
{\it Vertex operator algebras and the monster}, Pure and Applied Math, Vol.134,
Academic Press, London-New York, 1988.

\item{38.} {I.B.\ Frenkel} and {Y.\ Zhu}, ``Vertex
operator algebras associated to representations of affine and Virasoro
algebras'', {\it Duke Math.\ J.} {\bf 66} (1992), 123--168.

\item{39.} {J.\ Fuchs}, ``Fusion rules in conformal field theory'',
{\it Fortsch.\ Phys.} {\bf 42} (1994), 1--48.

\item{40.} {J.\ Fuchs, B.\ Gato-Rivera, B.\ Schellekens} and
{C.\ Schweigert}, ``Modular invariants and fusion rule
automorphisms from Galois theory'', {\it Phys.\ Lett.} {\bf B334} (1994), 113--120.

\item{41.} {J.\ Fuchs, A.N.\ Schellekens}, and {C.\ Schweigert},
``From Dynkin diagram symmetries to fixed point structures'', {\it Commun.\ Math.\
Phys.} {\bf 180} (1996), 39--97.

\item{42.} {J.\ Fuchs} and {C.\ Schweigert}, ``Branes:
from free fields to general backgrounds'', {\it Nucl.\ Phys.}
{\bf B530} (1998), 99--136.

\item{43.} L.\ Funar, ``On the TQFT representations of the mapping
class groups'', {\it Pac.\ J.\ Math.} {\bf 188} (1999), 251--274.

\item{44.} {P.\ Furlan, A.\ Ganchev,} and {V.\ Petkova}, 
``Quantum groups and fusion rules multiplicitites'', {\it Nucl.\ Phys.} {\bf
B343} (1990), 205--227.

\item{45.} {P.\ Furlan, A.\ Ganchev}, and {V.B.\ Petkova}, 
``An extension of the character ring of sl(3) and its
quantisation'', {\it Commun.\ Math.\ Phys.} {\bf 202} (1999), 701--733.

\item{46.} M.\ R.\ Gaberdiel, ``Fusion rules and logarithmic
representations of a WZW model at fractional level'', (hep-th/0105046).

\item{47.} {M.R.\ Gaberdiel} and {P.\ Goddard}, ``An
introduction to meromorphic conformal field theory and its
representations'', In: {\it Conformal Field Theory}, ed.\ by Y.\ Nutka et al,
Perseus Publishing, Cambridge MA, 2000.

\item{48.} M.R.\ Gaberdiel and H.G.\ Kausch, ``A rational logarithmic
conformal field theory'', {\it Phys.\ Lett.} {\bf B386} (1996), 131--137.

\item{49.} {T.\ Gannon}, ``Symmetries of the Kac-Peterson modular 
matrices of affine algebras'', {\it Invent.\ math.} {\bf 122} (1995), 341--357;

\item{} {T.\ Gannon, Ph.\ Ruelle,} and {M.\ A.\ Walton},
``Automorphism modular invariants of current algebras'', {\it Commun.\ Math.\ 
Phys.} {\bf 179} (1996), 121--156;

\item{} {T.\ Gannon}, ``The ADE$_7$ modular invariants of the
affine algebras'' (in preparation).

\item{50.} {T.\ Gannon}, ``The level 2 and 3 modular invariants
for the orthogonal algebras'', {\it Canad.\ J.\ Math.} {\bf 52} (2000),
503--521.

\item{51.} {T.\ Gannon}, ``The Cappelli-Itzykson-Zuber A-D-E
classification'', {\it Rev.\ Math.\ Phys.} {\bf 12} (2000), 739--748.

\item{52.} {T.\ Gannon}, ``The theory of fusion rings and modular
data'' (in preparation).

\item{53.} {T.\ Gannon}, ``Boundary conformal field theory and
fusion ring representations'', to appear in {\it Nucl.\ Phys.} {\bf B} (hep-th/016105).

\item{54.} {T.\ Gannon}, ``The modular invariants for the exceptional
algebras'' (in preparation);

\item{} {T.\ Gannon}, ``The modular invariants for the classical
algebras'' (work in progress).

\item{55.} F.\ R.\ Gantmacher, {\it The Theory of Matrices}, Chesea Publishing,
New York, 1990.

\item{56.} {G.\ Georgiev} and {O.\ Mathieu}, ''Cat\'egorie de
fusion pour les groupes de Chevalley'', {\it C.\ R.\ Acad.\ Sci.\ Paris} {\bf 315}
 (1992), 659--662.

\item{57.} {D.\ Gepner} and {E.\ Witten}, ``Strings on group
manifolds'', {\it Nucl.\ Phys.} {\bf B278} (1986), 493--549.

\item{58.} {F.\ M.\ Goodman} and {H.\ Wenzl},
``Littlewood-Richardson coefficients for Hecke algebras at roots of unity'',
 {\it Adv.\ Math.} {\bf 82} (1990), 244--265.

\item{59.} {A.\ Hanaki} and {I.\ Miyamoto},
  ``Classification of primitive association schemes of order up to
  22'', {\it Kyushu J.\ Math.}\ {\bf 54} (2000), 81--86.

\item{60.} {M.\ Hazewinkel, W.\ Hesselink, D.\ Siersma} and {
F.\ D.\ Veldkamp}, ``The ubiquity of Coxeter-Dynkin diagrams (an
introduction to the A-D-E problem)'', {\it Nieuw Arch.\ Wisk.} {\bf 25} (1977), 257--307;

\item{} {P.\ Slodowy}, 
``Platonic solids, Kleinian singularities, and Lie groups'', 
in: {\it Algebraic Geometry}, Lecture Notes
in Math 1008, Springer, Berlin, 1983, pp.102--138.

\item{61.} I.\ M.\ Isaacs, {\it Character Theory of Finite Groups}, Academic
Press, New York, 1976.

\item{62.} {V.\ G.\ Kac,} 
 ``Simple Lie groups and the Legendre symbol'', in: {\it Algebra}, 
 Lecture Notes in Math 848, Springer, New York, 1981, pp. 110--123.

\item{63.} {V.\ G.\ Kac,} {\it Infinite Dimensional Lie Algebras},
3rd edition, Cambridge University Press, Cambridge, 1990.

\item{64.} {V.G.\ Kac}, {\it Vertex Algebras for Beginners}, 2nd edn
(AMS, Providence, 1998).

\item{65.} {V.G.\ Kac} and {M.\ Wakimoto}, ``Modular invariant
representations of infinite-dimensional Lie algebras and superalgebras'',
{\it Proc.\ Natl.\ Acad.\ Sci.\ USA} {\bf 85} (1988), 4956--4960.

\item{66.} S.\ Kass, R.\ V.\ Moody, J.\ Patera, and R.\ Slansky, {\it
Affine Lie Algebras, Weight Multiplicities, and Branching Rules}, Vol.\ 1,
Univ.\ Calif.\ Press, Berkeley, 1990.

\item{67.} {C.\ Kassel}, {\it Quantum Groups}, Springer, New
York, 1995.

\item{68.} {T.\ Kawai}, ``On the structure of fusion
  algebras'', {\it Phys.\ Lett.}\ {\bf B217} (1989), 247--251.

\item{69.} {H.\ Kosaki, A.\ Munemasa} and {S.\ Yamagami},
``On fusion algebras associated to finite group actions'', {\it Pac.\ J.\ Math.}
{\bf 177} (1997), 269--290.

\item{70.} D.\ C.\ Lewellen, ``Sewing constraints for conformal field
theories on surfaces with boundaries'', {\it Nucl.\ Phys.} {\bf B372}
(1992), 654--682.

\item{71.} {G.\ Lusztig}, ``Exotic Fourier Transform'', {\it Duke
Math.\ J.} {\bf 73} (1994), 227--241.

\item{72.} {R.V.\ Moody} and {J.\ Patera},
``Computation of character decompositions of class functions on
campact semisimple Lie groups'', {\it Math.\ Comput.} {\bf 48} (1987),
799--827.

\item{73.} {G.\ Moore} and {N.\ Seiberg}, ``Classical and
quantum conformal field theory'', {\it Commun.\ Math.\ Phys.} {\bf 123}
(1989), 177--254.

\item{74.} {A.\ Munemasa}, ``Splitting fields of association schemes'',
{\it J.\ Combin.\ Th.} {\bf A57} (1991), 157--161.

\item{75.} {A.\ Ocneanu}, ``Paths on Coxeter diagrams: From Platonic
solids and singularities to minimal models and subfactors'', In: {\it Lectures 
on Operator Theory}. Amer.\ Math.\ Soc., Providence, 1999.

\item{76.}  {A.\ Ocneanu}, {\it (talks in Toronto and Kyoto, June
and December 2000)}.

\item{77.} {V.B.\ Petkova} and {J.-B.\ Zuber}, ``The
many faces of Ocneanu cells'', {\it Nucl.\ Phys.} {\bf B603} (2001), 449--496.

\item{78.} {K.-H.\ Rehren}, ``Braid group statistics and their
superselection rules'', In: {\it The Algebraic Theory of Superselection
Sectors}, D.\ Kastler (ed), World Scientific, Singapore, 1990, pp.333--355.

\item{79.} {P.\ Ruelle, E.\ Thiran,} and {J.\ Weyers},
``Implications of an arithmetic symmetry of the commutant for modular
invariants'', {\it Nucl.\ Phys.} {\bf B402} (1993), 693--708;

\item{} {M.\ Bauer, A.\ Coste, C.\ Itzykson} and {P.\
Ruelle}, ``Comments on the links between SU(3) modular invariants,
simple factors in the Jacobian of Fermat curves, and rational triangular
billiards'', {\it J.\ Geom.\ Phys.} {\bf 22} (1997), 134--189.

\item{80.} A.\ Sagnotti, Y.S.\ Stanev, ``Open descendents in conformal
field theory'', preprint (hep-th/9605042).

\item{81.} {A.\ N.\ Schellekens} and {S.\ Yankielowicz},
``Modular invariants from simple currents. An explicit proof.'' {\it
Phys.\ Lett.} {\bf B227} (1989), 387--391.

\item{82.} {A.\ Schilling} and {M.\ Shimozono},
``Bosonic formula for level restricted paths'', preprint (math.QA/9812106).

\item{83.} G.\ Segal, ``Geometric aspects of quantum field theories'', In:
Proc.\ Intern.\ Congr.\ Math., Kyoto, Springer, Hong Kong, 1991, pp.1387--1396.

\item{84.} {G.\ Tudose}, ``A special case of $sl(n)$-fusion
  coefficients'', preprint (math.CO/0008034).

\item{85.} {V.\ G.\ Turaev}, {\it Quantum Invariants of Knots and
3-Manifolds}, de Gruyter Studies in Math.\ 18, Berlin, 1994.

\item{86.} C.\ Vafa, {\it Phys.\ Lett.} {\bf 206B} (1988), 421--426.

\item{87.} {E.\ Verlinde}, ``Fusion rules and modular
  transformations in 2D conformal field theory'', {\it Nucl.\ Phys.}\
  {\bf B300} (1988), 360--376.

\item{88.} {M.\ A.\ Walton}, ``Algorithm for WZW fusion rules: a 
proof'', {\it Phys.\ Lett.} {\bf B241} (1990), 365--368.

\item{89.} {A.\ J.\ Wassermann}, ``Operator algebras and conformal
field theory'', In: {\it Proc.\ Intern.\ Congr.\ Math, Zurich},
Birkh\"auser, Basel, 1995, pp.966--979.

\item{90.} {N.\ J.\ Wildberger}, ``Duality and entropy of finite
commutative hypergroups and fusion rule algebras'', {\it J.\ London Math.\ Soc.}
{\bf 56} (1997), 275--291.

\item{91.} {E.\ Witten}, ``The Verlinde formula and the cohomology of
the Grassmannian'', In: {\it Geometry, Topology and Physics}, Conf.\ Proc.\
Lecture Notes in Geom.\ Top.\ vol.\ VI (1995), 357--422;

\item{} {B.\ Bertram, I.\ Ciocan-Fontanine} and {W.\
Fulton}, ``Quantum multiplication of Schur polynomials'', {\it J.\
Alg.} {\bf 219} (1999), 728--746.

\item{92.} {Y.\ Zhu}, ``Modular invariance of characters
of vertex operator algebras'', {\it J.\ Amer.\ Math.\ Soc.} {\bf 9}
(1996), 237--302.

\item{93.} {J.-B.\ Zuber}, ``Conformal, integrable and topological
theories, graphs and Coxeter groups'', In: {\it Proc.\ Intern.\ Conf.\
of Math.\ Phys.}, D.\ Iagnolnitzer (ed.), Intern.\ Press, Cambridge
MA, 1995, pp.674--689.

\bigskip\bigskip\bigskip\bigskip 
\centerline {{\bf Table.} Some affine
algebra modular invariant classifications}
$$\vbox{\tabskip=0pt\offinterlineskip
  \def\tablerule{\noalign{\hrule}}
 \halign to 4.5in{
    \strut#&\vrule#\tabskip=0em plus1em &    
    \hfil#&\vrule#&\hfil#&\vrule#&\hfil#&\vrule#&    
    \hfil#&\vrule#\tabskip=0pt\cr\tablerule      
&&algebra&&
 \omit\hidewidth \# of series \hidewidth&&  \omit\hidewidth levels of exceptionals\hidewidth
&&\omit\hidewidth verified for:\hidewidth &\cr\tablerule
&& $A_1^{(1)}$ && \hfil $k$ odd: 1 \hfil &&\hfil $k=10,16,28$\hfil &&\hfil$\forall k$\hfil& \cr
&& && \hfil $k$ even: 2 \hfil &&\hfil \hfil&& & \cr\tablerule
&& $A_2^{(1)}$ && \hfil $k$ arbitrary: 4 \hfil &&\hfil $k=5,9,21$\hfil&&\hfil$\forall k$\hfil & \cr\tablerule
&& $C_2^{(1)}$ && \hfil $k$ arbitrary: 2\hfil  &&\hfil $k=3,7,8,12$\hfil&&\hfil$k\le 25\,000$\hfil & \cr\tablerule
&& $G_2^{(1)}$ && \hfil $k$ arbitrary: 1 \hfil &&\hfil $k=3,4$\hfil&&\hfil$k\le 30\,000$\hfil & \cr\tablerule
&& $A_3^{(1)}$ && \hfil $k$ odd: 2\hfil  &&\hfil $k=4,6,8$\hfil&&\hfil$k\le 4000$\hfil & \cr
&&  && \hfil $k$ even: 4 \hfil &&&&& \cr\tablerule
&& $B_3^{(1)}$ && \hfil $k$ arbitrary: 2\hfil  &&\hfil $k=5,8,9$\hfil&&\hfil$k\le 3000$\hfil & \cr\tablerule
&& $C_3^{(1)}$ && \hfil $k$ odd: 1 \hfil &&\hfil $k=2,4,5$\hfil &&\hfil$k\le 4500$\hfil& \cr
&& && \hfil $k$ even: 2 \hfil &&\hfil \hfil&& & \cr\tablerule
&& $F_4^{(1)}$ && \hfil $k$ arbitrary: 1\hfil  &&\hfil $k=3,6,9$\hfil&&\hfil$k\le 2000$\hfil & \cr\tablerule
&& $E_6^{(1)}$ && \hfil $k$ arbitrary: 4\hfil  &&\hfil $k=4,6,12$\hfil&&\hfil$k\le 500$\hfil & \cr\tablerule
&& $E_7^{(1)}$ && \hfil $k$ odd: 1\hfil&&\hfil $k=3,12,18,(42?)$&&\hfil$k\le 400$\hfil &\cr
&& &&\hfil$k$ even: 2\hfil &&&&& \cr\tablerule
&& $E_8^{(1)}$ && \hfil $k$ arbitrary: 1\hfil  &&\hfil $k=4,30,(90?)$\hfil&&\hfil$k\le 380$\hfil & \cr
\tablerule\noalign{\smallskip}
}} $$

\end